\documentclass[reqno, 11pt]{amsart}
\usepackage{amscd, amssymb}

\theoremstyle{plain}
\newtheorem{Thm}[subsection]{Theorem}
\newtheorem{Cor}[subsection]{Corollary}
\newtheorem{Lem}[subsection]{Lemma}
\newtheorem{Prop}[subsection]{Proposition}
\newtheorem{Conj}[subsection]{Conjecture}

\theoremstyle{definition}
\newtheorem{Def}[subsection]{Definition}

\theoremstyle{remark}

\newtheorem{Rem}[subsection]{Remark}

\newtheorem{conj}{Conjecture}

\numberwithin{equation}{section}
\renewcommand{\rm}{\normalshape}

\newif\ifShowLabels
\ShowLabelstrue
\newdimen\theight
\def\TeXref#1{%
    \leavevmode\vadjust{\setbox0=\hbox{{\tt
        \quad\quad  {\small \rm #1}}}%
    \theight=\ht0
    \advance\theight by \lineskip
    \kern -\theight \vbox to
    \theight{\rightline{\rlap{\box0}}%
    \vss}%
    }}%

\ShowLabelsfalse

\renewcommand{\sec}[2]{\section{#2}\label{S:#1}%
    \ifShowLabels \TeXref{{S:#1}} \fi}
\newcommand{\ssec}[2]{\subsection{#2}\label{SS:#1}%
    \ifShowLabels \TeXref{{SS:#1}} \fi}

\newcommand{\refs}[1]{Section ~\ref{S:#1}}
\newcommand{\refss}[1]{Section ~\ref{SS:#1}}

\newcommand{\reft}[1]{Theorem ~\ref{T:#1}}
\newcommand{\refl}[1]{Lemma ~\ref{L:#1}}
\newcommand{\refp}[1]{Proposition ~\ref{P:#1}}

\newcommand{\refe}[1]{\eqref{E:#1}}

\newenvironment{thm}[1]%
    { \begin{Thm} \label{T:#1}  \ifShowLabels \TeXref{T:#1} \fi }%
    { \end{Thm} }

\renewcommand{\th}[1]{\begin{thm}{#1} \sl }
\renewcommand{\eth}{\end{thm} }

\newenvironment{lemma}[1]%
    { \begin{Lem} \label{L:#1}  \ifShowLabels \TeXref{L:#1} \fi }%
    { \end{Lem} }
\newcommand{\lem}[1]{\begin{lemma}{#1} \sl}
\newcommand{\elem}{\end{lemma}}

\newenvironment{propos}[1]%
    { \begin{Prop} \label{P:#1}  \ifShowLabels \TeXref{P:#1} \fi }%
    { \end{Prop} }
\newcommand{\prop}[1]{\begin{propos}{#1}\sl }
\newcommand{\eprop}{\end{propos}}

\newenvironment{corol}[1]%
    { \begin{Cor} \label{C:#1}  \ifShowLabels \TeXref{C:#1} \fi }%
    { \end{Cor} }
\newcommand{\cor}[1]{\begin{corol}{#1} \sl }
\newcommand{\ecor}{\end{corol}}

\newenvironment{defeni}[1]%
    { \begin{Def} \label{D:#1}  \ifShowLabels \TeXref{D:#1} \fi }%
    { \end{Def} }
\newcommand{\defe}[1]{\begin{defeni}{#1} \sl }
\newcommand{\edefe}{\end{defeni}}

\newenvironment{remark}[1]%
    { \begin{Rem} \label{R:#1}  \ifShowLabels \TeXref{R:#1} \fi }%
    { \end{Rem} }
\newcommand{\rem}[1]{\begin{remark}{#1}}
\newcommand{\erem}{\end{remark}}

\newenvironment{conjec}[1]%
    { \begin{Conj} \label{Co:#1}  \ifShowLabels \TeXref{Co:#1} \fi }%
    { \end{Conj} }
\renewcommand{\conj}[1]{\begin{conjec}{#1} \sl }
\newcommand{\econj}{\end{conjec}}

\newcommand{\eq}[1]%
    { \ifShowLabels \TeXref{E:#1} \fi
       \begin{equation} \label{E:#1} }
\newcommand{\eeq}{ \end{equation} }

\newcommand{\prf}{ \begin{proof} }
\newcommand{\epr}{ \end{proof} }

\newcommand\alp{\alpha}     

\newcommand\gam{\gamma}     \newcommand\Gam{\Gamma}
\newcommand\del{\delta}     \newcommand\Del{\Delta}

\newcommand\kap{\kappa}
\newcommand\lam{\lambda}        \newcommand\Lam{\Lambda}
\newcommand\sig{\sigma}     \newcommand\Sig{\Sigma}

     \newcommand\Ome{\Omega}

\newcommand\calD{{\mathcal{D}}}
\newcommand\calE{{\mathcal{E}}}
\newcommand\calF{{\mathcal{F}}}
\newcommand\calG{{\mathcal{G}}}
\newcommand\calH{{\mathcal{H}}}

\newcommand\calJ{{\mathcal{J}}}
\newcommand\calK{{\mathcal{K}}}
\newcommand\calL{{\mathcal{L}}}
\newcommand\calM{{\mathcal{M}}}

\newcommand\calO{{\mathcal{O}}}

\newcommand\calQ{{\mathcal{Q}}}
\newcommand\calR{{\mathcal{R}}}
\newcommand\calS{{\mathcal{S}}}
\newcommand\calT{{\mathcal{T}}}

\newcommand\calW{{\mathcal{W}}}
\newcommand\calX{{\mathcal{X}}}
\newcommand\calY{{\mathcal{Y}}}

\newcommand\PP{\mathbb{P}}
\renewcommand\AA{\mathbb{A}}

\newcommand\GG{\mathbb{G}}

\newcommand\ZZ{\mathbb{Z}}

\newcommand\CC{\mathbb{C}}

 \newcommand\grg{{\mathfrak{g}}}
 \newcommand\grh{{\mathfrak{h}}}

 \newcommand\grm{{\mathfrak{m}}}
 \newcommand\grn{{\mathfrak{n}}}
 
 \newcommand\grp{{\mathfrak{p}}}
 \newcommand\grq{{\mathfrak{q}}}
 \newcommand\grr{{\mathfrak{r}}}
 
 \newcommand\grt{{\mathfrak{t}}}

\newcommand\sdp{\times \hskip -0.3em {\raise 0.3ex
\hbox{$\scriptscriptstyle |$}}} 

\newcommand\Gr{\operatorname{Gr}}

\newcommand\Hom{\operatorname {Hom}}

\newcommand\id{\operatorname{id}}

\newcommand\im{\operatorname {im}}

\newcommand\rk{\operatorname{rk}}

\newcommand\Spec{\operatorname{Spec}}

\newcommand\Tr{\operatorname{Tr}}

\newcommand\hatT{{\widehat{T}}}

\newcommand\tilG{{\widetilde{G}}}

\newcommand\tilS{{\widetilde{S}}}

\newcommand\tilDel{{\widetilde{\Delta}}}

\newcommand\tilGam{{\widetilde{\Gamma}}}

\newcommand\tilpi{{\widetilde{\pi}}}

\newcommand\tilOme{{\widetilde{\Ome}}}

\newcommand\x{\times}
\newcommand\ten{\otimes}

\newcommand{\la}{\langle}

\newcommand\ch{\text{ch}}

\renewcommand\Spec{\operatorname{Spec}}

\newcommand\nc{\newcommand}
\nc{\coker}{\operatorname{coker}} \nc\tcF{\widetilde{\calF}}
\nc{\add}{\operatorname{add}}
\nc{\length}{\operatorname{length}}

\nc\hatLam{\widehat{\Lam}}

\nc{\Ups}{\Upsilon}
\nc{\tilUps}{\widetilde{\Ups}}
\nc{\kk}{\textsf{k}}
\nc\tgrp{\widetilde{\grp}}
\nc\tgrq{\widetilde{\grq}}

\nc\aff{\operatorname{aff}}
\nc\oGr{\overline{\Gr}} \nc\Bun{\operatorname{Bun}}
\nc\hgrg{\widehat{\grg}} \nc\Gal{\operatorname{Gal}}
\nc\tBl{\widetilde{\Bl}}
\nc\tilzeta{\widetilde{\zeta}}

\nc\tilBun{\widetilde \Bun} \nc\Pic{\operatorname{Pic}}
\nc\Rep{\operatorname{Rep}}
\nc\Bl{\operatorname{Bl}}
\nc{\val}{\operatorname{val}}
\nc{\tgrg}{\widetilde \grg}
\begin{document}
\title{The spherical Hecke algebra for affine Kac-Moody groups I}
\author{Alexander Braverman and David Kazhdan}
\begin{abstract}We define the spherical Hecke algebra for an (untwisted) affine Kac-Moody group over a local
non-archimedian field.
We prove a generalization of the Satake isomorphism for these algebras, relating it to integrable representations
of the Langlands dual affine Kac-Moody group. In the next publication we shall use these results to define
and study the notion of \emph{Hecke eigenfunction} for the group $G_{\aff}$.
\end{abstract}
\address{A.~B.~: Department of Mathematics, Brown University, 151 Thayer st., Providence RI, 02912, USA}
\address{D.~K.~: Department of Mathematics, the Hebrew University of Jerusalem, Givat Ram, Jerusalem, Israel}
\email{braval@math.brown.edu, kazhdan@math.huji.ac.il}
\maketitle
\sec{int}{Introduction}
\ssec{}{Langlands duality and the Satake isomorphism}
Let $F$ be a global field and let $\AA_F$ denote its ring of adeles.
Let $G$ be a split reductive group over $F$. The classical {\em Langlands duality} predicts that irreducible automorphic
representations of $G(\AA_F)$ are closely related to the homomorphisms from the absolute Galois group $\Gal_F$ of $F$ to
the {\em Langlands dual group} $G^{\vee}$. Similarly, if $G$ is a split reductive group over a local-nonarchimedian
field $\calK$, Langlands duality predicts a relation between irredicible representations of $G(\calK)$ and homomorphisms
from $\Gal_{\calK}$ to $G^{\vee}$.

The starting point for Langlands duality, which allows one to relate the "simplest" irreducible representations
of $G(\calK)$ (the so called {\em spherical representations}) to the "simplest" homomorphisms from $\Gal_{\calK}$ to $G^{\vee}$
(the unramified homomorphisms),
is the {\em Satake isomoprhism}, whose formulation we now recall.
Let $\calO\subset\calK$ denote
the ring of integers of $\calK$. Then the group $G(\calK)$ is a locally compact topological group and
$G(\calO)$ is its maximal compact subgroup. One may study the {\em spherical Hecke algebra} $\calH$ of
$G(\calO)$-biinvariant compactly supported $\CC$-valued measures on $G(\calK)$. The Satake isomorphism is a canonical
isomorphism between $\calH$ and the complexified
Grothendieck ring $K_0(\Rep(G^{\vee}))$
of finite-dimensional representations of $G^{\vee}$.
\ssec{}{The group $G_{\aff}$}To a connected reductive group $G$ as above one can associate the corresponding
affine Kac-Moody group $G_{\aff}$ in the following way.

Let $\Lam$ denote the coweight lattice of $G$ let $Q$ be an integral, even,
negative-definite symmetric bilinear form on $\Lam$
which is invariant under the Weyl group of $G$.

One can consider the polynomial loop group $G[t,t^{-1}]$ (this is an infinite-dimensional group ind-scheme).

It is well-known (cf. \cite{Kum}) that a form $Q$ as above gives rise to a central extension $\tilG$ of $G[t,t^{-1}]$:
$$
1\to \GG_m\to \tilG\to G[t,t^{-1}]\to 1
$$
(we are going to review the construction in \refss{central-ext}).
Moreover,\ $\tilG$ has again a natural structure of a group ind-scheme.

The multiplicative group $\GG_m$ acts naturally on $G[t,t^{-1}]$ and this action lifts to $\tilG$.
We denote the corresponding semi-direct product by $G_{\aff}$; we also let $\grg_{\aff}$ denote its Lie algebra.
Thus if $G$ is semi-simple then $\grg_{\aff}$ is an untwisted affine Kac-Moody Lie algebra in the sense of \cite{Kac};
in particular,
it can be described by the corresponding affine root system.

\ssec{}{The Hecke algebra for affine Kac-Moody groups}
Our dream is to develop some sort of Langlands theory in the case when $G$ is replaced by an affine Kac-Moody group
$G_{\aff}$. To do this we define the spherical Hecke algebra of $G_{\aff}$ in the following way.
Let $\calK$ and $\calO$ be as above. Then one may consider the group $G_{\aff}(\calK)$ and its subgroup
$G_{\aff}(\calO)$.

The group $G_{\aff}$ by definition maps to $\GG_m$; thus $G_{\aff}(\calK)$
maps to $\calK^*$. We denote this homomorphism by $\tilzeta$. In addition, the group $\calK^*$ is endowed with a natural
(valuation) homomoprhism to $\ZZ$. We denote its composition with $\tilzeta$ by $\pi$.

We now define the semigroup $G_{\aff}^+(\calK)$ to be the subsemigroup of $G_{\aff}(\calK)$ denerated by:

$\bullet$ the central $\calK^*\subset G_{\aff}(\calK)$;

$\bullet$ the subgroup $G_{\aff}(\calO)$;

$\bullet$ All elements $g\in G_{\aff}(\calK)$ such that $\pi(g)>0$.

We show (cf. \reft{main}(1)) that the convolution of any two double cosets of $G_{\aff}(\calO)$ inside
$G_{\aff}^+(\calK)$ is well-defined in the appropriate sense and in fact there is an associative algebra
structure on a suitable space of $G_{\aff}(\calO)$-biinvariant functions on $G^+_{\aff}(\calK)$. This algebra turns out
to be commutative and we call it {\it the spherical Hecke algebra of $G_{\aff}$} and denote it by $\calH(G_{\aff})$.
The algebra $\calH(G_{\aff})$ is graded by non-negative integers
(the grading comes from the map $\pi$ which is well-defined on double cosets
with resepct to $G_{\aff}(\calO)$); it is also an algebra over the field $\CC((v))$ of Laurent power series in a variable
$v$, which comes from the central $\calK^*$  in $G_{\aff}(\calK)$.
\ssec{}{The Satake isomorphism}The statement of the Satake isomorphism for $G_{\aff}$ is very similar to that for $G$.
First of all, in \refss{completed-invariant} we are going to define an analog of the algebra $\CC(T^{\vee})^W$ which we shall
denote by $\CC(\hatT^{\vee}_{\aff})^{W_{\aff}}$ (here $T_{\aff}=\CC^*\x T^{\vee}\x \CC^*$ is the dual of the maximal
torus of $G_{\aff}$, $W_{\aff}$ is the corresponding affine Weyl group and $\CC(\hatT^{\vee})$ denotes certain completion
of the algebra of regular functions on $T^{\vee}_{\aff}$). This is a finitely generated $\ZZ_{\geq 0}$-graded commutative algebra
over the field $\CC((v))$ of Laurent formal power series in the variable $v$ which should be thought of as a coordinate on
the third factor in $T^{\vee}_{\aff}=\CC^*\x T^{\vee}\x\CC^*$ (the grading has to do with the first factor); moreover, each component of the
grading is finite-dimensional over $\CC((v))$.

Assume that either $G$ is simply connected or that $G$ is a torus
\footnote{We believe that this assumption is not necessary, but at the moment we can't remove it}.
In this case we define (in \refs{main}) the {\it Langlands dual group} $G_{\aff}^{\vee}$. This is a group ind-scheme over $\CC$.
If $G$ is semi-simple, then
$G_{\aff}^{\vee}$ is another Kac-Moody group whose Lie algebra $\grg_{\aff}^{\vee}$ is an affine
Kac-Moody algebra with root system dual to that of $\grg_{\aff}$ (thus, in particular, it might be a twisted affine Lie algebra).
The group $G_{\aff}^{\vee}$ contains the torus $T_{\aff}^{\vee}$; moreover the first $\CC^*$-factor in $T^{\vee}_{\aff}$ is central
in $G^{\vee}_{\aff}$; also the projection $T^{\vee}_{\aff}\to\CC^*$ to the last factor extends to a homomorphism
$G_{\aff}^{\vee}\to \CC^*$.
It makes sense to consider integrable highest weight representations of $G_{\aff}^{\vee}$ and one can define
certain category $\Rep(G_{\aff}^{\vee})$ of such representations which is stable under tensor product (this category
contains all highest weight integrable representations of finite length, but certain infinite direct sums must
be included there as well). The complexified Grothendieck ring $K_0(G_{\aff}^{\vee})$ of this category is naturally isomorphic
to the algebra $\CC(\hatT^{\vee}_{\aff})^{W_{\aff}}$ via the character map. The corresponding grading
on $K_0(G_{\aff})$ comes from the central charge of $G_{\aff}^{\vee}$-modules and the action of the variable $v$
comes from tensoring $G_{\aff}^{\vee}$-modules by the one-dimensional representation coming from the homomorphism
$G_{\aff}^{\vee}\to\CC^*$, mentioned above.

The Satake isomorphism (cf. \reft{main}(2)) claims that the Hecke algebra $\calH(G_{\aff})$
is canonically
isomorphic to $\CC(\hatT^{\vee}_{\aff})^{W_{\aff}}$ (and thus when it makes sense also to $K_0(G_{\aff}^{\vee})$).

In the case when $G$ is semi-simple and simply connected, the group $G_{\aff}$ is an affine Kac-Moody group. We expect that with
slight modifications our Satake isomorphism should make sense for any symmetrizable Kac-Moody group. However, our proofs
are really designed for the affine case and do not seem to generalize to more general Kac-Moody groups.
\ssec{}{Further questions}It would definitely be interesting to extend the results of this paper to subgroups of
$G_{\aff}(\calK)$ other than $G_{\aff}(\calO)$. Arguably the most interesting of these subgroups is the analog of
the Iwahori subgroup of
$G_{\aff}(\calO)$. In this case we expect that the Hecke algebra is well-defined and will be closely related to
Cherednik's double affine
Hecke algebra. We plan to address this issue in \cite{BKM}.

Another natural question is this: by the definition, the algebra $\calH(G_{\aff})$ is endowed with a natural
basis given by the characteristic functions of
$G_{\aff}(\calO)\x G_{\aff}(\calO)$ orbits on $G_{\aff}(\calK)$. It would be quite interesting
to determine the image of this basis under the Satake isomorphism. When we deal with $G$ instead of $G_{\aff}$,
the answer is given by the so called {\em Hall polynomials} and it is essentially equivalent to the
"Macdonald formula for the spherical function" (cf. \cite{mac}).
A generalization of Macdonald formula to the affine case (closely related to the results  of \cite{CherMa} on affine
Hall polynomials) will
also be discussed in \cite{BKM}.

Let us also point out that the idea to interpret the convolution in $\calH(G_{\aff})$ in geometric terms (which is used
in order to prove the main results of this paper) came to us from the joint works of the first-named author with
M.~Finkelberg which give a partial generalization of the "geometric Satake correspondence" (cf. \cite{MV} and references therein)
to the affine case.
\ssec{}{Organization of the paper}This paper is organized as follows: in \refs{general} we introduce some notions
related to general Hecke algebras, in \refs{ex} we consider an example of (generalized) Hecke algebras (which later
turns out to be closely related to the algebra $\calH_{\aff}$ in the case when $G$ is a torus), in \refs{main}
we formulate our main \reft{main} describing the Hecke algebra $\calH_{\aff}$ in the general case.
The proof of \reft{main} occupies the last four sections of the paper. Although the statement of \reft{main}
is quite elementary, the proof uses heavily the machinery of algebraic geometry related to moduli spaces
of $G$-bundles on various algebraic surfaces (in particular, we don't know how to generalize our proofs to the
case when $G_{\aff}$ is replaced by a more general (non-affine) Kac-Moody group).
In addition,  in \refs{principal} we construct
explicitly the Satake isomorphism for $\calH_{\aff}$ by looking at its action on the principal series for
$G_{\aff}$ (in the spirit of \cite{Kap-H}).
\ssec{}{Acknowledgments}This work grew out of an attempt to find a "classical counterpart" of the
first author's joint work \cite{bf} with M.~Finkelberg whom we would like to thank for very useful
discussions and ideas; this paper has in fact grown out from an attempt to produce an elementary analog
of \cite{bf}. Also we are deeply grateful to V.~Drinfeld, who gave us many interesting ideas
on the subject. In addition we would like to thank I.~Cherednik, P.~Etingof, H.~Garland, V.~Kac, M.~Kapranov and
C.~Teleman for very illuminating
discussions related to the contents of this paper.

Part of this work has been carried out when the first author was visiting the
Hebrew University of Jerusalem. The work of both authors was partially supported by the BSF grant 5828239.
The first author was also partially supported by the NSF grant DMS-0600851.
\sec{general}{Generalities on Hecke algebras}
\ssec{good pairs}{Good pairs}
Let $\Gam$ be a group and let $\Gam_0\subset \Gam$ be a subgroup. Consider the multiplication map
$$
m:\Gam\underset{\Gam_0} \x \Gam\to \Gam.
$$
Let $X$ (resp. $Y$) be a subset of $\Gam$ which is right (resp. left) invariant
with respect to $\Gam_0$. Then we denote by $m_{X,Y}$ the restriction of $m$ to
$X\underset{\Gam_0}\x Y$.

We say that the pair $(\Gam,\Gam_0)$ is \emph{good} if for any two double cosets $X$ and $Y$
in $\Gam$ with respect to $\Gam_0$
we have:

1) The image of $m_{X,Y}$ consists of finite union of double cosets.

2) The map $m_{X,Y}$ has finite fibers.

\noindent
In this case one can define the Hecke algebra $\calH(\Gam,\Gam_0)$ as the convolution algebra of
$\Gam_0$-bi-invariant functions on $\Gam$ supported on finitely many double cosets with respect to $\Gam_0$.

\lem{good-equiv}
The pair $(\Gam,\Gam_0)$ is good if and only if for any $x\in \Gam$ the set
$\Gam_0x\Gam_0/\Gam_0$ is finite. Also, condition 2) above implies condition 1).
\elem
\prf
Assume that the condition of \refl{good-equiv} holds; let us prove that $(\Gam,\Gam_0)$ is a good pair.
Take any $x,y\in\Gam$ and set $X=\Gam_0x\Gam_0,Y=\Gam_0y\Gam_0$. Choose some $x_1,...,x_n\in X, y_1,...,y_k\in Y$
such that
$$
X=\bigsqcup\limits_{i=1}^n x_i\Gam_0,\quad Y=\bigsqcup\limits_{j=1}^k y_j\Gam_0
$$
Then $X\underset{\Gam_0}\x Y=\bigsqcup\limits_{i=1}^n x_i Y$. Since the restriction of $m_{X,Y}$ to each $x_iY$ is obviously injective
it follows that every fiber of $m_{X,Y}$ has at most $n$ elements. This proves condition 2).
On the other hand, $x_iY$ is the union of all the $x_iy_j\Gam_0$ and thus it follows that the image of $m_{X,Y}$ is
a finite union of right $\Gam_0$-cosets. In particular, it is a finite union of double cosets with respet to $\Gam_0$.

Let us now show that condition 2) implies that $\Gam_0x\Gam_0/\Gam_0$ is finite for all $x\in \Gam$ (in view of the
preceeding paragraph this will also imply that 2) implies 1)). Indeed, let $X=\Gam_0x\Gam_0$ and let
$Y=\Gam_0x^{-1}\Gam_0$. Consider the fiber of $m_{X,Y}$ over the unit element $e\in\Gam$. It is clearly equal
to $X/\Gam_0$ and thus condition 2) implies that it has to be finite.
\epr
The main source of example of good pairs is the following: assume that $\Gam$ is actually a locally compact
topological group and let $\Gam_0$ be any open compact subgroup (in particular, $\Gam$ can be a discrete
group and $\Gam_0$ - a finite subgroup).
\ssec{almost}{Generalization: almost good pairs}
We want to generalize the above notion in two directions. First, let $\Gam$ and $\Gam_0$ be as above and let
$\Gam^+\subset \Gam$ be a sub-semi-group containing $\Gam_0$. Then we may speak about the pair $(\Gam^+,\Gam_0)$ being
good and the Hecke algebra $\calH(\Gam^+,\Gam_0)$ makes sense.

\smallskip
\noindent
{\bf Remark.} In this case condition 2) no longer implies condition 1), since in the proof of this implication given above
we used inversion in $\Gam$.

\smallskip
\noindent
Let us also consider the following situation. Assume that we are given a central extension
$$
1\to A\to \tilGam\to \Gam\to 1
$$
of $\Gam$ by means of an abelian group $A$. Let also $B\subset A$ be a subgroup of $A$ such that
$A/B\simeq \ZZ$. Assume that $\Gam_0\subset \Gam$ can be lifted to a subgroup of $\tilGam_0\subset\tilGam$
fitting into the short exact sequence
$$
1\to B\to \tilGam_0\to \Gam_0\to 1.
$$
Note that in this case the group $\ZZ=A/B$ acts naturally on the set of double cosets of $\tilGam$ with respect to
$\tilGam_0$. We shall say that a subset $W\subset \tilGam$ is \emph{almost a finite union of double cosets
with respect to $\tilGam_0$} if $W$ is contained in the set $\ZZ_+(Z)$ where $Z\subset \Gam$ is a finite
union of double cosets.

Let  $\Gam^+\subset \Gam$ be a semi-group as above and $\tilGam^+$ is its preimage in $\tilGam$.
Then we say that the pair $(\tilGam^+,\tilGam_0)$ is \emph{almost good} if for any two double cosets
$X,Y$ of $\tilGam^+$
with respect to $\tilGam_0$ the  condition 2) above is satisfied
and in addition the following generalization of condition 1) is satisfied:

1') The image of the map $m_{X,Y}$ is almost a finite union of double cosets
with respect to $\tilGam_0$.

\noindent
{\bf Warning.} It is by no means true that if the pair $(\tilGam^+,\tilGam_0)$ is almost good then the pair
$(\Gam^+,\Gam_0)$ is good!

\smallskip
\noindent
In the above situation one may define the algebra $\calH(\tilGam^+,\tilGam_0)$ as the convolution
algebra of $\tilGam_0$-bi-invariant functions on $\tilGam^+$ whose support is almost a finite
union of double cosets. Note that this algebra can naturally be regarded as an algebra over formal
Laurent power series $\CC((v))$; here multiplication by $v$ corresponds to shifting the support of a function
by $1\in \ZZ$.

The algebra $\calH(\tilGam^+,\tilGam_0)$ has a natural subspace $\calH^{fin}(\tilGam^+,\tilGam_0)$ consisting of functions
supported on finitely many double cosets; it is naturally a module over the ring $\CC[v,v^{-1}]$ of Laurent
polynomials. By the definition, we have
$$
\calH(\tilGam^+,\tilGam_0)=\calH^{fin}(\tilGam^+,\tilGam_0)\underset{\CC[v,v^{-1}]}\ten \CC((v)).
$$
In general $\calH^{fin}(\tilGam^+,\tilGam_0)$ is not a subalgebra
of $\calH(\tilGam^+,\tilGam_0)$.
Let $R$ be any subring of $\CC((v))$ containing $\CC[v,v^{-1}]$. We say that
$\calH(\tilGam^+,\tilGam_0)$ is \emph{defined over $R$} if $\calH^{fin}(\tilGam^+,\tilGam_0)\underset{\CC[v,v^{-1}]}\ten R$
is a subalgebra of $\calH(\tilGam^+,\tilGam_0)$. In this case we set
$\calH_R(\tilGam^+,\tilGam_0)=\calH^{fin}(\tilGam^+,\tilGam_0)\underset{\CC[v,v^{-1}]}\ten R$.

\ssec{}{Good actions}
Let $(\Gam,\Gam_0)$ be a (not necessarily good) pair and let $\Ome$ be a $\Gam$-set. Let
$a: \Gam\underset{\Gam_0}\x\Ome\to \Ome$ be the natural map, coming from the action of $\Gam$ on $\Ome$.
Let now $X$ as before be a subset of $\Gam$ which is right-invariant under $\Gam_0$ and let $Z$ be a
$\Gam_0$-invariant subset of $\Ome$. Then we shall denote by $a_{X,Z}$ the restriction of $a$ to
$X\underset{\Gam_0}\x Z$.

 We say that the action of $\Gam$ on $\Ome$
is good if for any double coset $X$ in $\Gam$ with respect to $\Gam_0$ and for any $\Gam_0$-orbit
$Z$ in $\Ome$ the following conditions are satisfied:

$1_{\Ome}$: The image of $a_{X,Z}$ is a union of finitely many $\Gam_0$-orbits.

$2_{\Ome}$: The map $a_{X,Z}$ has finite fibers.

\noindent
Let $\calF(\Ome)$ denote the space of $\Gam_0$-invariant $\CC$-valued functions on $\Ome$, which are supported
on finitely many $\Gam_0$-orbits (in other words, $\calF(\Ome)$ is the space of functions on
$\Gam_0\backslash\Ome$
with finite support). Then conditions $1_{\Ome}$ and $2_{\Ome}$ above guarantee that we can define an action
of $\calH(\Gam,\Gam_0)$ on $\calF(\Ome)$ by putting
$$
h(f)(z)=\sum\limits_{(g,w)\in\Gam\underset{\Gam_0}\x\Ome, ~ a(g,w)=z} h(g)f(w).
$$
Here $h\in\calH(\Gam,\Gam_0)$, $f\in \calF(\Ome)$, $z\in\Ome$.

The following lemma is straightforward.
\lem{}
Let $(\Gam,\Gam_0)$ be any pair and let $\Ome$ be any $\Gam$-set. Assume that condition $2_{\Ome}$ is satisfied.
Then the pair $(\Gam,\Gam_0)$ satisfies condition 2 from \refss{good pairs}.
\elem
\ssec{}{Almost good actions}Let $(\tilGam,\tilGam^+,\tilGam_0)$ be as in \refss{almost} and let $\Ome$ as before
be a $\Gam$-set. Let also $\tilOme$ be an $A$-torsor over $\Ome$ on which $\tilGam$ acts in such a way that
$A\subset \tilGam$ acts on it in the natural way and the action of $\tilGam$ on $\tilOme$ is compatible with
the action of $\Gam$ on $\Ome$. Then we may define the notion of an almost good action of $\tilGam^+$ on $\tilOme$
(similarly to the notion of an almost good pair).
Namely, we say that a $\tilGam_0$-invariant subset $W$ of $\tilOme$ is almost a finite union of $\tilGam_0$-orbits
(or just {\it almost finite} for brevity) if
$W$ is contained in $\ZZ_+(Z)$ where $Z$ is a finite union of $\tilGam_0$-orbits (note that $\ZZ=A/B$ acts naturally
on the set of $\tilGam_0$-orbits in $\tilOme$). We say that the action of $\tilGam$ on $\tilOme$ is almost
good if condition $2_{\tilOme}$ is satisfied together with the following condition $1'_{\tilOme}$:

\smallskip
$1'_{\tilOme}$: For any $\tilGam_0$ double coset $X$ in $\tilGam^+$ and for any $\tilGam_0$-orbit $Z$ in $\tilOme$,
the image of $a_{X,Z}$ is almost a finite union of $\tilGam_0$-orbits.

\smallskip
\noindent
We now let $\calF^{fin}(\tilOme)$ denote the space of $\tilGam_0$-invariant $\CC$-valued functions on $\tilOme$
which are supported on a finite union of $\tilGam_0$-orbits and we also let
$\calF(\tilOme)$ denote the space of $\tilGam_0$-invariant $\CC$-valued functions on $\tilOme$ which are supported on
an almost finite union of $\tilGam_0$-orbits. It is easy to see that if the action of $\tilGam^+$ on $\tilOme$
is almost good, then the algebra $\calH(\tilGam^+,\tilGam_0)$ acts on $\calF(\tilOme)$.
\ssec{goodactgen}{A generalization}Let $(\Ome,\tilOme)$ be as above and let $(\Ome',\tilOme')$ be another such pair.
Assume that we are given an $\tilGam$-equivariant morphism $\varpi:\tilOme\to\tilOme'$. We shall say that a
$\tilGam_0$-invariant subset $Z$ of $\tilOme$ is {\it almost finite with respect to $\tilOme'$} if

a) $\varpi(Z)$ is almost finite;

b) The natural map $\tilGam_0\backslash Z\to\tilGam_0\backslash\varpi(Z)$ has finite fibers.

\noindent
We shall say that $\tilOme$ is {\it almost good with respect to $\tilOme'$} if the following two conditions
hold:

(i) For any $X\in \tilGam_0\backslash \tilGam^+/\tilGam_0$ and any $Z\in\tilGam_0\backslash \tilOme$ the image of the map
$a_{X,Z}$ is almost finite with respect to $\tilOme'$.

(ii) $\tilOme'$ is almost good.

\noindent
We claim that condition (ii) implies the condition $2_{\tilOme}$. Indeed, without loss of generality we may that
$\Ome$ and $\Ome'$ are homogenous spaces for $\Gam$ (resp. $\tilOme$ and $\tilOme'$ are homogenous spaces for
$\tilGam$), i.e. $\Ome=\Gam/\Del,\tilOme=\tilGam/\tilDel,\Ome'=\Gam/\Del',\tilOme=\tilGam/\tilDel'$, where
$\Del\subset\Del'\subset \Gam$ and $\tilDel\subset\tilDel'\subset \tilGam$. Let
$X$ be a double coset in $\tilGam^+$ with respect to $\tilGam_0$. Choose
$z,\gam\in \tilGam$. Set $X=\tilGam_0x\tilGam_0$, $Z=\tilGam_0 z\tilDel/\tilDel$ and
$Z'=\tilGam_0 z\tilDel/\tilDel'$. Then the fiber of $a_{X,Z}$ over the image of $\gam$ in $\tilOme$ consists of
all triples $(x,\del)\in X\x\tilDel$ such that $xz\del=\gam$. Similarly,
the fiber of $a_{X,Z'}$ over the image of $\gam$ in $\tilOme'$ consists of
all triples $(x,\del')\in X\x\tilDel'$ such that $xz\del'=\gam$. Since
$\tilDel\subset\tilDel'$ it follows that the former fiber is embedded in the latter. Thus, the fact that
the fibers of $a_{X,Z'}$ are finite implies that the fibers of $a_{X,Z'}$ are finite.

In addition it is obvious that condition (ii) implies part a) of condition (i).
Let now $\calF_{\tilOme'}(\tilOme)$ denote
the space of all $\Gam_0$-invariant functions on $\tilOme$ whose support is almost finite with respect to $\tilOme'$.
Then it is easy to see that the conditions (i) and (ii) imply that
the Hecke algebra $\calH(\tilGam^+,\tilGam_0)$ acts on $\calF_{\tilOme'}(\tilOme)$.
\sec{ex}{An example}
\ssec{}{The setup}Let $\Lam$ be a lattice of finite rank (i.e. $\Lam$ is an abelian group isomorphic to $\ZZ^l$ for
some $l$). Let $\Gam= (\Lam\oplus\Lam)\rtimes \ZZ$ where $\ZZ$ acts on $\Lam\oplus\Lam$ by means
of the autormorphism
$\sig$ sending $(\lam,\mu)$ to $(\lam,\lam+\mu)$. We shall usually write an element
of $\Gam$ as $(\lam,\mu,k)$. Define the subgroup $\Gam_0\subset\Gam$ as the subgroup consisting of all elements
of the form $(\lam,0,0)$.

Define also the semigroup $\Gam^+$ by setting
$$
\Gam^+=\{ (\lam,0,0)|\ \lam\in \Lam\}~\bigcup~
\{(\lam,\mu,k)\ \text{where $\lam$ and $\mu$ are arbitrary and $k>0$}\}.
$$

Let $b(\cdot,\cdot)$ be any $\ZZ$-valued bilinear form on $\Lam$. We set
$$
Q(\lam,\mu)=b(\lam,\mu)+b(\mu,\lam).
$$
Then $Q$ is a symmetric $\ZZ$-valued form on $\Lam$. Moreover, $Q$ is even, i.e. $Q(\lam,\lam)$
is even for any $\lam\in\Lam$. We define the central extension $\widehat{\Lam\oplus \Lam}$ of
$\Lam\oplus\Lam$ by means of $\ZZ$ in the following way: $\widehat{\Lam\oplus \Lam}=\ZZ\x(\Lam\oplus\Lam)$
as a set and the multiplication is defined by
$$
(a_1,\lam_1,\mu_1)(a_2,\lam_2,\mu_2)=
(a_1+a_2+b(\lam_1,\mu_2)-b(\mu_1,\lam_2),\lam_1+\lam_2,\mu_1+\mu_2).
$$
It is well known that $\widehat{\Lam\oplus \Lam}$ depends (canonically) only on $Q$ and not on $b$.

The automorphism $\sig$ of $\Lam\oplus\Lam$ considered above extends naturally to
$\widehat{\Lam\oplus \Lam}$. Abusing the notation, we shall denote the resulting
automorphism of $\widehat{\Lam\oplus \Lam}$ also by $\sig$. It is defined by
$$
\sig:(a,\lam,\mu)\mapsto (a,\lam,\lam+\mu).
$$
We define $\tilGam=\widehat{\Lam\oplus \Lam}\rtimes \ZZ$ where the generator of
$\ZZ$ acts on $\widehat{\Lam\oplus \Lam}$ by means of $\sig$. Explicitly,
we have $\tilGam=\ZZ\x(\Lam\oplus\Lam)\x\ZZ$ as a set and the multiplication is given by
\eq{multip}
\begin{aligned}
&(a_1,\lam_1,\mu_1,k_1)(a_2,\lam_2,\mu_2,k_2)=\\
&=(a_1+a_2+b(\lam_1,\mu_2)-b(\mu_1,\lam_2)+k_1b(\lam_1,\lam_2),\lam_1+\lam_2,\mu_1+\mu_2+k_1\lam_2,k_1+k_2).
\end{aligned}
\eeq
We set $\tilGam_0$ to consist of all elements of $\tilGam$ of the form $(0,\lam,0,0)$. Note that in this
case $\tilGam_0$ is naturally isomorphic to $\Gam_0$.

As before we define $\tilGam^+$ to be the preimage of $\Gam^+$ in $\tilGam$.
\ssec{}{The variety $X_{Q,v}$}From now on we assume that $Q$ is non-degenerate.
Let $\Lam^{\vee}$ denote the dual lattice to $\Lam$. Then $Q$ defines a homomorphism $e:\Lam\to\Lam^{\vee}$,
which is injective since  $Q$ is non-degenerate.
Let $T=\Lam\ten\CC^*,~T^{\vee}=\Lam^{\vee}\ten\CC^*$.
Then $\Lam^{\vee}$ is a lattice of characters of $T$ and
$\Lam$ is the lattice of characters of $T^{\vee}$.

For every $v\in\CC^*$ and any $\lam^{\vee}\in\Lam^{\vee}$ we denote by $v^{\lam^{\vee}}$ the corresponding
element of $T^{\vee}$. Assume that $|v|<1$. In this case we set
$$
X_{Q,v}=T^{\vee}/v^{e(\Lam)}.
$$
It is well-known that $X_{Q,v}$ is a complex abelian variety.
The form $Q$ also defines a holomorphic line bundle $\calL_{Q,v}$ on $X_{Q,v}$ in the following way.
Let $\pi:T^{\vee}\to X_{Q,v}$ be the natural projection.
Then for any open subset $U$ of $X_{Q,v}$ we define the space of holomorphic sections of
$\calL_{Q,v}$ to be the space of holomorphic functions $F$ on $\pi^{-1}(U)$ satisfying
$$
F(xv^{e(\nu)})=v^{-\frac{Q(\nu,\nu)}{2}}x^{-\nu}F(x).
$$
Here $x\in T^{\vee}$ and $x^{-k\nu}=\nu(x)^k$.
The line bundle $\calL_{Q,v}$ is ample if and only if $Q$ is negative definite.
In this case we set
$$
\CC[X_{Q,v}]=\bigoplus\limits_{k=0}^{\infty}\Gam(X_{Q,v},\calL_{Q,v}^{\ten k}).
$$
Then $\CC[X_{Q,v}]$ is a  graded holomorphic vector bundle of algebras on the punctured unit disc
$\calD^*=\{v\in\CC^*,\ |v|<1\}$. More precisely, let $R$ denote the subring of $\CC((v))$ consisting of those formal series
which are convergent for $0<|v|<1$. Then there exists a pair $(X_Q,\calL_Q)$ where
$X_Q$ is an abelian variety over $R$ and $\calL_Q$ is a line bundle over it, such that the specialization
of $(X_Q,\calL_Q)$ to any $v\in\calD^*$ is isomorphic to $(X_{Q,v},\calL_{Q,v})$. We set $\CC[X_Q]$ to be the
direct image of the sheaf to $Spec(R)$; it can be naturally regarded as a graded $R$-algebra.
We also set $\widehat{\CC[X_Q]}=\CC[X_Q]\underset{R}\ten\CC((v))$.

Below is the main result of this Section.

\th{ex-main}
Assume that $Q$ is negative-definite. Then
\begin{enumerate}
\item
The pair $(\tilGam^+,\tilGam_0)$ defined above is almost good.
\item
There is a natural isomorphism
$$
\calH(\tilGam^+,\tilGam_0)\simeq \widehat{\CC[X_Q]}
$$
of graded $\CC((v))$-algebras.
\item
The algebra $\calH(\tilGam^+,\tilGam_0)$ is defined over the ring $R$ and we have the natural
isomorphism
$$
\calH_R(\tilGam^+,\tilGam_0)\simeq \CC[X_Q]
$$
of graded $R$-algebras.
\end{enumerate}
\eth
\prf
The proof is a straightforward calculation. Let us first describe $\calH(\tilGam^+,\tilGam_0)$
as a vector space. It follows from \refe{multip} that we have the following formula:
\eq{multip'}
(0,\nu,0,0)(a,\lam,\mu,k)=(a+b(\nu,\mu),\lam+\nu,\mu,k).
\eeq
Hence it follows that every left $\tilGam_0$-coset contains unique element of the form
$(a,0,\mu,k)$. Let us investigate when two such elements lie in the same double coset
with respect to $\tilGam_0$. For any $\nu\in\Lam$ we have
$$
(a,0,\mu,k)(0,\nu,0,0)=(a-b(\mu,\nu),\nu,\mu+k\nu,k).
$$
Multiplying the result on the left by $(0,-\nu,0,0)$ and applying \refe{multip'} again
we get
$$
(a-b(\nu,\mu)-b(\mu,\nu)-kb(\nu,\nu),0,\mu+k\nu,k)=(a-Q(\mu,\nu)-\frac{kQ(\nu,\nu)}{2},0,\mu+k\nu,k).
$$
In other words, we see that an element of $\calH(\tilGam^+,\tilGam_0)_k$ is a
function $f$ on $\ZZ\x\Lam$ ("corresponding" to $a$ and $\mu$) satisfying the relation
\eq{relation}
f(a,\mu)=f(a-Q(\mu,\nu)-\frac{k Q(\nu,\nu)}{2},\mu+k\nu)
\eeq
for any $\nu\in\Lam$ and such that:

1) For fixed $a$ the number $f(a,\mu)$ is non-zero only for finitely many values of $\mu$.

2) There exists some $a_0\in\ZZ$ such that $f(a,\mu)=0$ for any $\mu$ and any $a<a_0$.
\footnote{Note that here we use the assumption that $Q$ is negative-definite}

Set now
$$
F(v,x)=\sum\limits_{(a,\mu)\in\ZZ\x\Lam} f(a,\mu)v^a x^{\mu}.
$$
Then \refe{relation} together with 1) and 2) above imply that $F$ can naturally be regarded
as an element of $\widehat{\CC[X_Q]}_k$. In other words, we get an isomorphism of graded
vector spaces $\calH(\tilGam^+,\tilGam_0)\simeq \widehat{\CC[X_Q]}$.
Let us prove that this is actually an isomorphism of algebras.

Let $L$ be the collection of all the elements of $\tilGam$ of the form $(a,0,\mu,k)$. It is easy
to see that $L$ is actually a normal subgroup of $\tilGam$ isomorphic to $\ZZ\x\Lam\x\ZZ$, which can also be identified with
the kernel of the homomorphism $\tilGam\to \Lam$ sending $(a,\lam,\mu,k)$ to $\lam$.
It is also easy to see that $L$ and $\Gam_0$ satisfy the following properties:

\medskip
1) $L\cap \Gam_0=\{ e\}$.

2) $\tilGam=L\cdot \Gam_0=\Gam_0\cdot L$.

\medskip
\noindent
In other words $L$ is a normal subgroup of $\tilGam$ whose elements provide unique representatives for both left and right
cosets of $\tilGam$ with respect to $\Gam_0$. In particular, since $L$ is normal, we have a natural (conjugation) action
of $H$ on $L$. Let $\CC[L]$ denote the space $\CC$-valued functions on $L$ with
almost finite support. This is an algebra with respect to convolution (this algebra is a completion of the group algebra of $L$). Let also
$\CC[L]^H$ be the space  $H$-invariants in $\CC[L]$ (note that it follows from \refe{relation} that every $H$-orbit in
$L$ is almost finite). Then it is easy to see that conditions 1) and 2) above imply that the restriction map from $\tilGam$ to
$L$ defines an isomorphism of algebras
$\calH(\tilGam^+,\tilGam_0)\simeq \CC[L]^H$.
However, $L$ is isomorphic to $\ZZ\x\Lam\x\ZZ$, which implies that $\CC[L]^H\simeq \widehat{\CC[X_Q]}$.
\epr
\ssec{reformulation}{A reformulation}The above result can be reformulated as follows.
 Define an action of $\Lam$ on the torus $T^{\vee}\x \CC^*$ by the following formula:
\eq{formula}
\nu(x,t)=(v^{e(\nu)}x, t x^{\nu} v^{\frac{Q(\nu,\nu)}{2}}).
\eeq
Then the quotient $(T^{\vee}\x\CC^*)/\Lam$ can be naturally identified with the total space of the complement to
the zero section in the line
$\calL_{Q,v}^{-1}$ on $X_{Q,v}$. In particular, the algebra $\CC[X_{Q,v}]$ is equal to the algebra of regular
functions on $(T^{\vee}\x\CC^*)/\Lam$. In other words, let us consider the torus
$T^{\vee}_{\aff}=\CC^*\x T^{\vee}\x \CC^*$ with "coordinates" $(t,x,v)$. There is a natural action of
$\Lam$ on $T^{\vee}_{\aff}$, defined by the obvious analog of \refe{formula}:
\eq{formula'}
\nu(t,x,v)=( t x^{\nu} v^{\frac{Q(\nu,\nu)}{2}},v^{e(\nu)}x,v).
\eeq

Also, for every $k\in\ZZ$ let $\CC(T^{\vee}_{\aff})_k$ denote the space of regular functions on $T^{\vee}_{\aff}$
which are homogeneous of degree $k$ with respect to the first $\CC^*$. Let
$\CC(\hatT^{\vee}_{\aff})_k$ denote the completion of this space in the $v$-adic topology. We also set
\eq{completed-pol}
\CC(\hatT^{\vee}_{\aff})=\bigoplus\limits_{k\in\ZZ}\CC(\hatT^{\vee}_{\aff})_k.
\end{equation}
 Then there is natural identification
$\widehat{\CC[X_Q]}\simeq \CC(\hatT^{\vee}_{\aff})^{\Lam}$ (it is easy to see that $\CC(T^{\vee}_{\aff})_k^{\Lam}\neq 0$
if and only if $k\geq 0$).
\sec{main}{The main result}
\ssec{}{Loop groups and their cousins}
In this paper for convenience we adopt the polynomial version of loop groups (as opposed to formal
loops version; cf. however, \refss{variant}, which explains how to reformulate the story in terms of formal
loops).

Let $G$ be a split connected reductive algebraic group over a field $\kk$ with Lie algebra $\grg$.
Let $T$ be a maximal (split) torus in $G$ and let $\grt$ denote the corresponding Cartan subalgebra
of $\grg$. We also denote by $W$ the Weyl group of $G$.


Let $\Lam$ denote the coweight lattice of $G$. Note that we can regard $\Lam$ as a subset of $\grh$.
We shall assume that there exists an integral, negative-definite symmetric bilinear form $Q$ on $\Lam$
which $W$-invariant. In this case $Q\ten \kk$ is a restriction to $\Lam$ of a $G$-invariant form on $\grg$,
which we shall denote by $Q_{\grg}$.

Consider the corresponding polynomial Lie algebra $\grg[t,t^{-1}]$ and the group $G[t,t^{-1}]$.
As was mentioned in the introduction, the form $Q$ gives rise to a central extension $\tilG$ of $G[t,t^{-1}]$:
$$
1\to \GG_m^*\to \tilG\to G[t,t^{-1}]\to 1.
$$
Moreover, both $G[t,t^{-1}]$ and $\tilG$ has a natural structure of a group ind-scheme over $\kk$
\footnote{This group-scheme is non-reduced if $G$ is not semi-simple};
We denote by $\hgrg$ the corresponding Lie algebra; it fits into the exact sequence
$$
0\to k\to {\widetilde \grg}\to\grg[t,t^{-1}]\to 0.
$$
The group $\GG_m$ acts naturally on $G[t,t^{-1}]$ and this action lifts to $\tilG$.
We denote the corresponding semi-direct product by $G_{\aff}$; we also let $\grg_{\aff}$ denote its Lie algebra.
Thus if $G$ is semi-simple then $\grg_{\aff}$ is an untwisted affine Kac-Moody Lie algebra in the sense of \cite{Kac};
in particular,
it can be described by the corresponding affine root system.

We also let $G_{\aff}'$ denote the semi-direct product $\GG_m\ltimes G[t,t^{-1}]$.
\ssec{completed-invariant}{The algebra $\CC(\hatT^{\vee}_{\aff})^{W_{\aff}}$}
Let us recall that for $G$ as above one can consider the Langlands dual group $G^{\vee}$. One of its crucial properties
says that the complexified Grothendieck group $K_0(\Rep G^{\vee})$ of the category $\Rep(G^{\vee})$ of finite-dimensional
representations is naturally isomorphic to the algebra $\CC(T^{\vee})^W$ of $W$-invariant polynomials
on the torus $T^{\vee}$.

We would like to define analogous notions when $G$ is replaced by $G_{\aff}$. First of all, let us describe the analog
of the algebra $\CC(T^{\vee})$. Similarly to \refs{ex} let us set $T_{\aff}=\GG_m\x T\x\GG_m$ and $T^{\vee}_{\aff}=
\CC^*\x T^{\vee}\x \CC^*$ (we shall only work with the dual torus $T^{\vee}_{\aff}$ over $\CC$). It will be convenient for
us to think about the first $\GG_m$-factor in $T_{\aff}$ as dual to the {\em second} $\CC^*$-factor in $T_{\aff}^{\vee}$ (and vice
versa).
We also define $\CC(\hatT^{\vee}_{\aff})$ as in \refe{completed-pol}.

The weight lattice $\Lam_{\aff}$ of $T^{\vee}_{\aff}$ is naturally identified with $\ZZ\x\Lam\x \ZZ$. It is also the coweight lattice
of $T_{\aff}$ (however, let us stress that according to our conventions the first multiple in $\Lam_{\aff}=\ZZ\x\Lam\x \ZZ$ corresponds
to the loop rotation in $G_{\aff}$ and the last multiple corresponds to the center of $G_{\aff}$).
Let $\Lam_{\aff}'=\ZZ\x\Lam$. For any $k\in \ZZ$ we denote by $\Lam_{\aff,k}$ (resp. by $\Lam_{\aff,k}'$) the set of all elements
of $\Lam_{\aff}$ (resp. of $\Lam_{\aff}'$ whose first coordinate is equal to $k$.

Let $W_{\aff}$ denote affine Weyl group of $G$ which is the semi-direct product of $W$ and $\Lam$.
It acts on the lattice $\Lam_{\aff}$ (resp. $\Lam_{\aff}'$)
preserving each $\Lam_{\aff,k}$ (resp. each $\Lam_{\aff,k}'$). In order to describe this action explicitly it is convenient
to set $W_{\aff,k}=W\ltimes k\Lam$ which naturally acts on $\Lam$.
 Then the restriction of the $W_{\aff}$-action to $\Lam_{\aff,k}\simeq\Lam\x\ZZ$
comes from the natural $W_{\aff,k}$-action on the first multiple.

In the case when $G$ is semi-simple and simply connected  the set $\Lam_{\aff,k}/W_{\aff}$ admits the following description.
Let us denote by $\Lam_{\aff}^+$ the set of dominant
coweights of $G_{\aff}$ (i.e. the set of those coweights whose scalar product with every positive
root of $G_{\aff}$ is non-negative) . We also denote by $\Lam_{\aff,k}$ be the subset of $\Lam_{\aff}$ consisting
of all elements of the form $(k,\lam,n)$. We put $\Lam_{\aff,k}^+=\Lam_{\aff}^+\cap \Lam_{\aff,k}$.

Let
$\Lam_k^+\subset \Lam$ denote the set of dominant coweights of $G$ such that $\la \lambda,\alp)\leq k$
when $\alp$ is the highest root of $\grg$.
Then it is well-known that a weight $(k,\lam,n)$ of
$G_{\aff}^{\vee}$ lies in $\Lam_{\aff,k}^+$  if and only if $\lam\in\Lam_k^+$ (thus $\Lam_{\aff,k}^+=\Lam_k^+\x \ZZ$).

It is well known that every $W_{\aff}$-orbit on $\Lam_{\aff,k}$ contains unique element of $\Lam_{\aff,k}^+$.
This is equivalent to saying that $\Lam_k^+\simeq \Lam/W_{\aff,k}$.

\ssec{}{The Langlands dual group}
We would like to interpret the algebra $\CC(\hatT^{\vee}_{\aff})^{W_{\aff}}$ as the Grothendieck ring of a certain category
representation of the appropriately defined affine Langlands dual group $G_{\aff}^{\vee}$. Unfortunately, we don't know a good
definition of $G_{\aff}^{\vee}$ for general $G$; however, we can define $G_{\aff}^{\vee}$ (and the appropriate category of representations) in the following
two cases:

1) $G$ is torus

2) $G$ is simply connected.

\noindent
Let us explain these definitions (we are not going to use the group $G_{\aff}^{\vee}$ in the remaining part of the paper; however, we think that
interpreting  the algebra $\CC(\hatT^{\vee}_{\aff})^{W_{\aff}}$ (that will one of the most important players in our main \reft{main}) at least in some
cases is very instructive).

\medskip
\noindent
{\em Case 1: $G$ is a torus.}

\noindent
Let us use the notations of the previous section.
In particular, we let $e:\Lam\to\Lam^{\vee}$ denote the map given by the form $Q$. We let $G^{\vee}_Q$ denote the
torus over $\CC$ whose lattice of cocharacters is $e(\Lam)$. Thus $G^{\vee}_Q$ fits into the short exact sequence
$$
0\to Z\to G^{\vee}_Q\to G^{\vee}\to 0,
$$
where $Z=\Lam^{\vee}/e(\Lam)$.
Note that the form $Q$ can be used to define the dual form $Q^{\vee}$ on $e(\Lam)$ which is also
integral and even. Thus we can form the group $(G^{\vee}_Q)_{\aff}$. It contains the group
$G^{\vee}_Q$ as a subgroup; it is easy to see that $Z\subset G^{\vee}_Q$ is central in
$(G^{\vee}_Q)_{\aff}$. Thus we may define $G^{\vee}_{\aff}=(G^{\vee}_Q)_{\aff}/Z$.

The quotient of $G_{\aff}^{\vee}$ by the central $\CC^*$ can be described
as follows. Consider the group ind-scheme $(G^{\vee})_{\aff}'$. It is easy to see that its connected components are numbered
by $\Lam^{\vee}$. Then $G_{\aff}^{\vee}/\CC^*$ is equal to the union of the above connected components corresponding to
the elements of $e(\Lam)\subset\Lam^{\vee}$.

\medskip
\noindent
{\em Case 2: $G$ is semi-simple and simply connected.}

\noindent
In this case let $\grg_{\aff}^{\vee}$ denote the {\it Langlands dual} affine Lie algebra of $\grg_{\aff}$ (considered as a Lie algebra
over $\CC$. By definition, this is an affine
Kac-Moody Lie algebra whose root system is dual to that of $\grg_{\aff}$. Moreover, it comes also with a fixed central
sublagebra
$\CC\subset\grg_{\aff}^{\vee}$ which is determined uniquely by the choice of $Q_{\grg}$. Note that in general (when $\grg$ is not
simply laced) the algebra $\grg_{\aff}^{\vee}$ is not isomorphic to $(\grg^{\vee})_{\aff}$ (here $\grg^{\vee}$ as before denotes
the Langlands dual Lie algebra of $\grg$). Moreover, if $\grg$ is not simply laced, then $\grg_{\aff}^{\vee}$ is a twisted
affine Lie algebra. However, the algebra $\grg_{\aff}^{\vee}$ always contains $\grg^{\vee}\x \CC^2$ as a Levi subalgebra.

\smallskip
\noindent
We now let $G_{\aff}^{\vee}$ denote any connected group ind-scheme over $\CC$ such that

a) The Lie algebra of $G_{\aff}^{\vee}$ is $\grg_{\aff}^{\vee}$

b) The above embedding $\grg^{\vee}\hookrightarrow \grg_{\aff}^{\vee}$ extends to an embedding
$G^{\vee}\hookrightarrow G_{\aff}^{\vee}$.

\smallskip
\noindent
The existence (and uniqueness) of $G_{\aff}^{\vee}$ follows immediately from \cite{Tits}.

In both cases 1) and 2) the group $G_{\aff}^{\vee}$ contains the torus $T_{\aff}^{\vee}=$
Here the first $\ZZ$-factor
is responsible for the center of $G_{\aff}^{\vee}$;
it can also be thought of as coming from the loop
rotation in $G_{\aff}$. The second $\ZZ$-factor is responsible for the loop rotation in $G_{\aff}^{\vee}$.
Similarly, we shall denote by $\Lam_{\aff}'=\Lam\x\ZZ$ the weight lattice of $G_{\aff}'$.

It follows easily from properties 1),2) and 3) that
the group $G_{\aff}^{\vee}$ maps naturally to $\CC^*$
(this homomorphism is dual to the central embedding $\CC^*\to
G_{\aff}$). We denote the kernel of this homomorphism by
$\tilG^{\vee}$ and we let $\tgrg^{\vee}$ denote its Lie algebra.
It is clear that in fact $G_{\aff}^{\vee}=\tilG^{\vee}\rtimes\CC^*$.
We shall fix such an isomorphism, which, in particular, endows $G_{\aff}^{\vee}$
with a subgroup $\CC^*$, which we shall call the loop rotation subgroup in $G_{\aff}^{\vee}$.

Also $G_{\aff}^{\vee}$ has a natural central subgroup isomorphic to $\CC^*$. In addition, it is easy to see from
condition 3) that $G_{\aff}^{\vee}$ is generated by $(G_{\aff}^{\vee})^0$ and by $T_{\aff}^{\vee}$.

\ssec{}{Representations of $G_{\aff}^{\vee}$}
The Lie algebra $\grg_{\aff}$ is endowed with a canonical "parabolic" subalgebra $\grp_0$ which is equal
to $\kk\oplus\grg[t]\oplus\kk$. We let $\grp_0^{\vee}$ denote the corresponding dual subalgebra of
$\grg_{\aff}^{\vee}$ and we set $\grn_0^{\vee}=[\grp_0^{\vee},\grp_0^{\vee}]$.

By an integrable highest weight representation of $G_{\aff}^{\vee}$ we shall mean an algebraic representation
of $G_{\aff}^{\vee}$ whose restriction to the Lie algebra $\grn_0^{\vee}$ is locally nilpotent.
It is easy to see that any such representation is semi-simple. We say that such a representation $L$ is of level $k$
if the central $\CC^*$ acts on $L$ by means of the character $z\mapsto z^k$. Then it is well-known that

a) If $L$ is a non-zero representation of level $k$ then $k\geq 0$.

b) Any representation of level $0$ is a pull-back of a representation of $\CC^*$ under the above-mentioned map
$G_{\aff}^{\vee}\to\CC^*$.

c) The irreducible representations of level $k>0$ are in one-to-one correspondence with elements of
$\Lam_{\aff,k}/W_{\aff}$.

Let $\Rep(G_{\aff}^{\vee})$ denote the category of integrable highest weight
representations $L$ of $G_{\aff}^{\vee}$ satisfying the following conditions:

\medskip
(i) $L$ is a direct sum of irreducible representations of $G_{\aff}^{\vee}$ with finite multiplicities.

(ii) For every $n\in\ZZ$ let us denote by $L_n$ the subspace of $L$ on which the loop rotation subgroup
$\CC^*$ acts by means of the character $z\mapsto z^n$. Then $L_n=0$ for $0\ll n$.

\medskip
\noindent
The category $\Rep(G_{\aff}^{\vee})$ is stable under tensor product (cf. \cite{Kac}). Also the group $\ZZ$ acts naturally
on this category by multiplying every $L$ by the corresponding character of $G_{\aff}^{\vee}$ coming from the homomorphism
$G_{\aff}^{\vee}\to \GG_m$. We denote by $K_0(G_{\aff}^{\vee})$ the complexified
Grothendieck ring of $\Rep(G_{\aff}^{\vee})$ tensored with $\CC$. It is clear that the above
$\ZZ$-action gives rise to a $\CC((v))$-module
structure on $K_0(G_{\aff}^{\vee})$. In addition the algebra $K_0(G_{\aff}^{\vee})$ is $\ZZ_+$-graded:
by the definition, elements of degree $k$ in $K_0(G_{\aff}^{\vee})$ correspond to
representations of $G_{\aff}^{\vee}$ of level $k$.
In this way, $K_0(G_{\aff}^{\vee})$ becomes a graded $\CC((v))$-algebra.

To every $L\in \Rep(G_{\aff}^{\vee})$ we can consider its character $\ch(L)$ which is an element of $\CC(\hatT^{\vee}_{\aff})$.
The assignment $L\mapsto \ch(L)$ extends to an isomorphism
\eq{k-polyn}
K_0(G_{\aff}^{\vee})\simeq \CC(\hatT^{\vee}_{\aff}.
\end{equation}
In the future we are going to work with the algebra $\CC(\hatT^{\vee}_{\aff})$, which is defined for any reductive $G$.
However, we think that intuitively it is important keep in mind the isomorphism \refe{k-polyn}.
\ssec{sphint}{The spherical Hecke algebra of $G_{\aff}$}
Let now $\kk$ be a finite field. We shall now choose a ring $\calK$ and its subring $\calO$, which will
come from one of the following two situations:

\medskip
a) $\calK$ is a local non-archimedian field with residue field $\kk$ and $\calO$ is its ring of integers

b)  $\calK=\kk[s,s^{-1}]$, $\calO=\kk[s]$.

\medskip
\noindent
Case a) is somewhat more interesting and natural. However, we still include case b), since our proofs
are a little more transparent
in this case (although the difference between the two cases is not essential). In both cases we are given the natural
valuation homomorphism $val:\calK^*\to\ZZ$, whose kernel is $\calO^*$.

It now makes sense to consider
the group $\tilGam=G_{\aff}(\calK)$ and its subgroup $\tilGam_0=G_{\aff}(\calO)$
(we also have $\Gam=G'_{\aff}(\calK),\Gam_0=G'_{\aff}(\calO)$.
Note that this
pair satisfies the conditions of \refss{almost} with $A=\calK^*$ which (in both cases a) and b)) has a natural map to
$\ZZ$ with kernel being $B=\calO^*$.

We have the natural homomorphisms
$$
\zeta:\Gam=G'_{\aff}(\calK)\to \calK^*\quad\text{and}\quad val:\calK^*\to\ZZ.
$$
We set $\pi=val\circ\zeta$. Similarly we have the homomorphisms
$\tilzeta:G_{\aff}(\calK)\to\calK^*$ and
$\tilpi:G_{\aff}(\calK)\to \ZZ$.

We set
$$
\Gam^+=G'_{\aff}(\calO)\cup \pi^{-1}(\ZZ_{>0}).
$$
As in \refss{almost} we let $\tilGam^+$ denote the preimage of $\Gam^+$ in $\tilGam$.

Below is the main result of this paper:
\th{main}
\begin{enumerate}
\item
The pair $(\tilGam^+,\tilGam_0)$ defined above is almost good in the terminology
of \refss{almost}.
\item
The Hecke algebra $\calH(\tilGam^+,\Gam_0)$ is naturally
isomorphic to $\CC(\hatT_{\aff}^{\vee})$ as a
graded $\CC((v))$-algebra.
\end{enumerate}
\eth
In fact, we are going to construct the isomorphism (2) explicitly;
we shall call it the Satake isomorphism for the group $G_{\aff}$. We expect that it should be possible
to describe this isomorphism explicitly in terms of representation theory of $G_{\aff}^{\vee}$ in the spirit
of \cite{Lu-qan} (cf. also \cite{bf} for some closely related conjectures). However we shall postpone the detailed
discussion of this question for another publication.

The proof of the first assertion occupies \refs{bundles}, \refs{determinant} and \refs{central}. Although the first
assertion seems to be quite elementary, the only proof of \reft{main}(1) that we know requires
interpreting the fibers of the corresponding maps $m_{X,Y}$ in terms of algebraic geometry.
The proof of the second assertion occupies
\refs{principal}.

\sec{bundles}{The algebra $\calH(\tilGam^+,\tilGam_0)$ via bundles on surfaces}
\ssec{}{Double cosets and Kleinian singularities}

For any $k>0$ let $S_k=\Spec \kk[x,y,z]/xy-z^k$ and let $S^0_k$ be the complement to the point
$(0,0,0)$ in $S_k$. This is a smooth surface.

Let $\mu_k$ denote the group-scheme of roots of unity of order $k$. Then the natural map
$\grp_k:\AA^2\to S_k$ given by
$$
(u,v)\mapsto (u^k,v^k,uv)
$$
identifies
$S_k$ with $(\AA^2\backslash \{ 0\})/\mu_k$ (note that the action of $\mu_k$ on $(\AA^2\backslash \{ 0\})$
is free).

Let $\Sig$ be some other surface over $\kk$ and let $a$ be a point of $\Sig$. We shall say that
$x$ is of type $A_k$ if there is an etale neighborhood of $a$ in $\Sig$ which is isomorphic to an etale
neighborhood of $(0,0,0)$ in $S_k$. In particular, if $k=1$ this will just mean that $a$ is a smooth point of
$\Sig$.

\ssec{}{A variant: the case of an arbitrary local field}
The surface $S_k$ introduced above will play a crucial role in the proof of \reft{main} in the case b) from \refss{sphint}.
Let us explain how to define it in the case a), i.e. when we work over  a local non-archimedian field $\calK$ with ring
of integers $\calO$. Let $z$ be a uniformizer in $\calO$ (it is convenient to choose it, though nothing
will actually depend on this choice). In this case we shall set $S_k=\Spec\calO[x,y]/xy-z^k$.This a scheme over $\calO$.
It is easy to see that if $k=1$ then this scheme is regular (warning: it is, however, NOT
regular as a scheme over $\calO$). We set $S_k^0$ to be the complement in $S_k$ of the (closed) point
given by the equations $x=y=z=0$. The scheme $S_k^0$ is regular for all $k>0$.
We define a map $\grp_k:S_1\to S_k$  of $\calO$-schemes which sends $(u,v)$ to
$(x=u^k,y=v^k)$. It is clear that under this identification we get
$S_k^0=S_1^0/\mu_k$ where $\mu_k$ is the group-scheme of $k$-th roots of unity (over $\calO$).

In all the proofs in the rest of this paper we shall always assume that we are in case b) (in the terminology
of \refss{sphint}). The extension to case a) will
be straightforward using the notations introduced in this subsection.

We shall say that a two-dimensional scheme $\Sig$ has a singularity of type $A_k$ at a closed point $p$ if near $p$ it is
etale-locally isomorphic to $S_k$ with $p$ corresponding to the point $(0,0,0)$ (in particular, if $k=1$ this means that
$p$ is a smooth point of $\Sig$).

\medskip
\noindent
For any $\kk$-variety $S$ we shall denote by $\Bun_G(S)$ the set of isomorphism classes of principal $G$-bundles
on $S$.
\prop{bijectionw}
The following sets are in natural bijection:
\begin{enumerate}
\item
The set $\Bun_G(S^0_k)$;
\item
The set $G'_{\aff}(\calO)\backslash\pi^{-1}(k)/G'_{\aff}(\calO)$;
\item
The set $\Lam/W_{\aff,k}$;
\item
The set of $G$-conjugacy classes of homomorphisms $\mu_k\to G$.
\end{enumerate}
\eprop
\prf
Let us first establish the bijection between (1) and (4). Let $\calF$ denote a $G$-bundle on $S_k$.
Consider the $G$-bundle $\grp_k^*(\calF)$ on $\AA^2\backslash\{ 0\}$. It extends uniquely to the whole
of $\AA^2$ and thus it is trivial. On the other hand, the bundle $\grp_k^*(\calF)$ is $\mu_k$-equivariant.
Since this bundle is trivial, such an equivariant structure gives rise to a homomorphism $\mu_k\to G$
defined uniquely up to conjugacy. This defines a map $\Bun_G(S^0_k)\to \Hom(\mu_k,G)/G$. It is clear that
this is actually a bijection, since a $\mu_k$-equivariant $G$-bundle on $\AA^2\backslash\{ 0\}$
descends uniquely to $S_k$.

Let us now establish the bijection between (1) and (2). To do that let us denote by $\sig$ the automorphism
of $G(\calK)$ sending $g(t,s)$ to $g(ts,s)$. Let us now identify $\pi^{-1}(k)$ with $G(\calK)$ with right
$G(\calK)$ action being the standard one (by right shifts) and with left $G(\calK)$-action given by
$g(h)=\sig^K_0(G)h$. Thus we have the bijection
\eq{maindoubleq}
G'_{\aff}(\calO)\backslash\pi^{-1}(k)/G'_{\aff}(\calO)=
G[ts^k,t^{-1}s^{-k},s]\backslash G[t,t^{-1},s,s^{-1}]/G[t,t^{-1},s].
\eeq
Let $U_k=\Spec \kk[ts^k,t^{-1}s^{-k},s]$, $V_k=\Spec \kk[t,t^{-1},s]$. Both $U_k$ and $V_k$ are isomorphic to
$\GG_m\x \AA^1$ and thus every $G$-bundle on either of these surfaces is trivial. Both $U_k$ and $V_k$ contain
$W=\Spec \kk[t,t^{-1},s,s^{-1}]\simeq \GG_m\x \GG_m$ as a Zariski open subset.
Thus the RHS of \refe{maindoubleq} can be identified with $\Bun_G(S')$ where $S'$ is obtained by gluing
$U_k$ and $V_k$ along $W$. Hence it remains to construct an isomorphism
$S'{\widetilde \to} S_k$. Such an isomorphism can be obtained by setting $x=ts^k, y=t^{-1}, z=s$.

It remains to construct a bijection between (3) and (4).
Recall that $\Lam=\Hom(\GG_m,T)$. Thus, since $\mu_k$ is a closed subscheme of $\GG_m$,
given every $\lam\in \Lam$ we may restrict it to $\mu_k$ and get a homomorphism $\mu_k\to T$.
By composing it with the embedding $T\hookrightarrow G$ we get a homomorphism $\mu_k\to G$ which clearly
depends only on the image of $\lam$ in $\Lam/W_{\aff,k}$. Thus we get a well-defined map
$\Lam/W_{\aff,k}\to (\Hom(\mu_k,G)/G)$. The surjectivity of this map follows from the fact that
$\mu_k$ is a diagonalizable group-scheme. For the injectivity note that since any two elements in $T$ which are
conjugate in $G$ lie in the same $W$-orbit in $T$. Thus it is enough to show that
for any two homomorphisms $\lam,\mu:\GG_m\to T$ whose restrictions to $\mu_k$ coincide the difference,
$\lam-\mu$ is divisible by $k$. This is enough to check for $T=\GG_m$ where it is obvious.
\epr
\ssec{}{Groupoids}Recall that a groupoid is a category in which all morphisms are isomorphisms.
All groupoids that will appear in this paper will be small. A typical example of a groupoid
for us will be like this: if a group
$H$ acts on a set $X$, then we can consider the groupoid $X/H$ where the objects are points of $X$ and a morphism
from $x$ to $y$ is an elements $h\in H$ such that $h(x)=y$.
Using the terminology common in the theory of algebraic stacks,
we shall say that a functor $f:\calY\to\calY'$ between two groupoids is representable
if any $y\in\calY$ has no non-trivial automorphisms which act trivially on $f(y)$. In this case
for any $y'\in \calY'$ we define the fiber $f^{-1}(y')$ to be the set of isomorphism classes
of pairs $(y,\alp)$, where $y\in\calS$ and $\alp$ is an isomorphism between $f(y)$ and $y'$.

By a groupoid structure on a set $Y$ we shall mean a groupoid $\calY$ whose set of isomorphism classes
is identified with $Y$.

The double quotient
$G'_{\aff}(\calO)\backslash\pi^{-1}(k)/G'_{\aff}\calO)$ and the set
$\Bun_G(S^0_k)$ have natural groupoid structures
and it is easy to see that the above identification between the them is actually an equivalence of groupoids.

From now on we shall treat these sets as groupods (abusing the notation, we shall denote them in the same
way as before).
\ssec{convolution}{The convolution diagram}
Choose any two positive integers $k$ and $l$.
Consider the following "convolution diagram":
\eq{conv}
G'_{\aff}\calO)\backslash\pi^{-1}(k)\underset{G'_{\aff}\calO)}\x \pi^{-1}(l)/G'_{\aff}\calO)\overset{m}\to
G'_{\aff}\calO)\backslash\pi^{-1}(s^{k+l})/G'_{\aff}\calO),
\eeq
where $m$ denotes the multiplication map. As was mentioned above we shall think about this diagram as
a map of groupoids. We want to give an interpretation of this diagram in geometric terms,
i.e. in terms of $G$-bundles on some surfaces. The RHS of \refe{conv} is clearly equal to
$\Bun_G(S^0_{k+l})$. Let us interpret the LHS.

The surface $S_{k+l}$ has canonical crepant resolution $\tilS_{k+l}$. Its fiber over $0\in S_{k+l}$ consists
of a tree of $k+l-1$ rational curves $E_1,...,E_{k+l-1}$ where the self-intersection of each $E_j$ is $-2$.
Thus it is possible to blow down all the $E_j$ except $E_k$. Let us call the resulting surface $S_{k,l}$.
It has two potentially singular points corresponding to $0$ and $\infty$ in $E_k$. They are of type $A_k$ and $A_l$
respectively (i.e. these points are really singular if the corresponding integer ($k$ or $l$) is greater than
1).
More precisely, we claim that $S_{k,l}$ can be covered by two open subsets which are isomorphic to $S_k$,
$S_l$ respectively. Let us describe how this is done in more detail.

The surface $S_k$ can be covered by two open subsets $U_k$ and $V_k$, each one isomorphic to
$\AA^1\x \GG_m$. Namely if $S_k=\Spec \kk[x,y,z]/xy-z^k$ then
$U_k=\Spec \kk[x,x^{-1},z]$ and $V_k=\Spec[y,y^{-1},z]$ (note that $U_k$ and $V_k$ can be considered as open subsets
both $S_k$ and $S_k^0$). Similarly, let us consider
$S_l=\Spec \kk[u,v,w]/uv-w^l$. Let us identify $V_k$ with $U_l$ by setting $z=w$ and $y=u^{-1}$.
Let us now glue $S_k$ and $S_l$ along $V_k$ and $U_l$ respectively. The resulting surface $\Sig$ maps
to $S_{k+l}=\Spec\kk[p,q,r]/pq-r^{k+l}$ in the following way: we map $S_k$ to $S_{k+l}$ by sending
$(x,y,z)$ to $(x,yz^l,z)$ and we map $S_l$ to $S_{k+l}$ by sending $(u,v,w)$ to $(uw^k,v,w)$.
These formulas are compatible, since by multiplying the equality $y=u^{-1}$ by
$xu$ on both sides we get $uxy=uz^k$ in the left hand side and $x$ in the right hand side, which implies
that $x=uz^k=uw^k$. Similarly, multiplying $u=y^{-1}$ by $vy$ on both sides, we get the equality
$yz^l=v$. These equalities make sense as long as $y$ and $u$ are invertible.
It is now easy to see that the resulting map $\Sig\to S_{k+l}$ is an isomorphism away from the point
$(0,0,0)$ and the fiber over $(0,0,0)$ is naturally isomorphic to $\PP^1$.
From this it is easy to deduce that $\Sig\simeq S_{k,l}$.
We let $\grp_{k,l}$ denote the natural map $S_{k,l}\to S_{k+l}$.
Note that $S_{k,l}\backslash \grp_{k,l}^{-1}(0,0,0)\simeq S_{k+l}$.

We let $S_{k,l}^0$ denote the complement to the two singular points in $S_{k,l}$.

\smallskip
\noindent

\prop{conv-no-cent}
a) The LHS of \refe{conv} can be naturally identified with $\Bun_G(S^0_{k,l})$ (as a groupoid);

b) The corresponding map $\Bun_G(S^0_{k,l})\to \Bun_G(S^0_{k+l})$ is just the restriction to the complement of
$\grp_{k,l}^{-1}(0,0,0)$.
\eprop
\prf
According to the proof of \refp{bijectionw} the convolution diagram \refe{conv}
can be identified with
\eq{conv'}
\begin{aligned}
&G[ts^{k+l},t^{-1}s^{-k-l},s]\backslash G[t,s,t^{-1},s^{-1}]\underset{G[ts^l,t^{-1}s^{-l},s]}\x
G[t,t^{-1},s,s^{-1}]/G[t,t^{-1},s]\overset{m}\to\\
&G[ts^{k+l},t^{-1}s^{-k-l},s]\backslash G[t,s,t^{-1},s^{-1}]
/G[t,t^{-1},s]
\end{aligned}
\eeq
We claim now, that the LHS of \refe{conv'} can be naturally identified with $\Bun_G(S^0_{k,l})$. Indeed, let us set
$$
x=ts^{k+l}, y=t^{-1}s^{-l}, u=ts^l, v=t^{-1}, z=w=s.
$$
Then we have $xy=z^k$ and $uv=w^l$ and
$\kk[ts^{k+l},t^{-1}s^{-k-l},s]=\kk[x,x^{-1},z]$, $\kk[ts^l,t^{-1}s^{-l},s]=\kk[y,y^{-1},z]$.
Thus the product
$$
G[ts^{k+l},t^{-1}s^{-k-l},s]\backslash G[t,t^{-1},s,s^{-1}]\x
G[t,t^{-1},s,s^{-1}]/G[t,t^{-1},s]
$$ classifies $(\calF_k,\calF_l,\beta_k,\alp_l)$
where

a) $\calF_k$ is a $G$-bundle on $S_k$ and $\beta_k$ is a trivialization of $\calF_k$ on $V_k$;

b) $\calF_l$ is a $G$-bundle on $S_l$ and $\beta_l$ is a trivialization of $\calF_l$ on $U_l$.

\noindent
It follows from this that the product
$$
G[ts^{k+l},t^{-1}s^{-k-l},s]\backslash G[t,s,t^{-1},s^{-1}]\underset{G[ts^l,t^{-1}s^{-l},s]}\x
G[t,t^{-1},s,s^{-1}]/G[t,t^{-1},s]
$$
classifies the triples $(\calF_k,\calF_l,\gam)$ where $\calF_k$ and $\calF_l$ are as above and $\gam$
is an isomorphism between $\calF_k|_{V_k}$ and $\calF_l|_{U_l}$ (with respect to the identification of
$V_k$ and $U_l$ discussed in the previous subsection). But such a triple is the same as a $G$-bundle
$\calF_{k,l}$ on $S_{k,l}$ since the latter surface is obtained by gluing $S_k$ and $S_l$ by
identifying $V_k$ with $U_l$.
\epr
\ssec{}{A generalization}Fix now some positive integers $k_1,...,k_n$ and let $k=k_1+\cdots k_n$.
In this case we can define the surface $S_{k_1,...,k_n}$ which is obtained by gluing
the surfaces $S_{k_1},\cdots ,S_{k_n}$ by identifying $V_{k_1}$ with $U_{k_2}$, $V_{k_2}$ with $U_{k_3}$ etc.
We have the map $\grp_{k_1,\cdots,k_n}:S_{k_1,\cdots,k_n}\to S_k$. The surface $S_{k_1,\cdots,k_n}$ is smooth
away from $n$ singular points which all lie in the preimage of zero under $p_{k_1,\cdots,k_n}$ and
are of types $A_{k_1},\cdots,A_{k_n}$ respectively. We shall denote by $S_{k_1,\cdots,k_n}^0$ the complement
to these singular points.

It is easy to see that away from $\grp_{k_1,\cdots,k_n}^{-1}(0,0,0)$ the map
$\grp_{k_1,\cdots,k_n}$ is an isomorphism. In particular, $S_k^0$ embeds into $S_{k_1,\cdots,k_n}^0$ as an open subset.

On the other hand, we can consider the map (of groupoids)
$$
\begin{aligned}
 G'_{\aff}(\calO)\backslash \pi^{-1}(k_1)\underset{G'_{\aff}(\calO)}\x\pi^{-1}(k_2)\underset{G'_{\aff}(\calO)}\x\cdots
\underset{G'_{\aff}(\calO)}\x \pi^{-1}(k_n)/G'_{\aff}(\calO)&\to \\
G'_{\aff}(\calO)\backslash \pi^{-1}(k)&/G'_{\aff}(\calO).
\end{aligned}
$$
Then (by using the same argument as above) we can identify this map with the restriction map
$$
\Bun_G(S^0_{k_1,\cdots,k_n})\to \Bun_G(S^0_k).
$$
\ssec{}{Fibers of the convolution morphism}
Consider again the convolution morphism \refe{conv}. It is easy to see that this morphism is representable and we
would like to understand its fibers. According to \refss{convolution}
we can instead look at the fibers of the restriction map
$\Bun_G(S^0_{k,l})\to \Bun_G(S^0_{k+l})$. By the definition, the fiber of this map over a bundle $\calF_{k+l}$ consists
of all extensions $\calF_{k,l}$ of $\calF_{k+l}$ to $S_{k,l}$. If the pair $(\Gam^+,\Gam_0)$ were good, this would
imply that this set is finite. However, it is easy to see that this set is actually infinite. In
\refs{central} we shall see that this problem can be remedied by changing the group $G'_{\aff}$ by $G_{\aff}$.
First we need to recall some (mostly well-known) constructions from algebraic geometry.
\sec{determinant}{Determinant torsors}
\ssec{determinants}{Relative determinants}
Let  $X$ be a smooth variety (over an arbitrary field) and let $\calQ$ be a vector bundle on $X$.
Then we shall denote by $\det(\calQ)\in\Pic(X)$ the top exterior power of $\calQ$. More generally, let
$\calQ$ be an arbitrary coherent sheaf on $X$. Then locally $\calQ$ has a resolution
$$
0\to\calE_n\to\cdots\to\calE_1\to\calE_0\to\calQ\to 0
$$
where all $\calE_i$ are locally free. In this case
 we set
$$
\det(\calQ)=\bigotimes\limits_{i=0}^n \det(\calE_i)^{(-1)^i}
$$
It is well-known that the result is canonically independent
of the choice of the resolution (thus, in particular, it makes sense globally).

Let now $S$ be a smooth variety  and let $X$ be a smooth divisor in $S$ (in the future
$S$ will usually be a surface and $X$ will usually be a closed curve inside $S$).
Set $S^0=S\backslash X$.

Let $\calF_1,\calF_2$ be two locally free sheaves on $S$ together with an isomorphism $\alp:\calF_1|_{S^0}
\widetilde{\to}\calF_2|_{S^0}$. In this case we can form the \emph{relative determinant}
$\det_X(\calF_1,\calF_2,\alp)\in\Pic(X)$
\footnote{In the sequel we shall just write $\det_X(\calF_1,\calF_2)$ or $\det_X(\alp)$ when it does not lead to a confusion}
of $\calF_1$ and $\calF_2$ in the following way.

Assume first, that $\alp$ defines an embedding of $\calF_1$ into
$\calF_2$ as coherent sheaves. Then we may consider the quotient $\calQ=\calF_2/\calF_1$. This is a coherent sheaf on $S$
which is set-theoretically concentrated on $X$. Thus we can find a filtration
$0\subset \calQ_1\subset \calQ_2\subset\cdots\subset \calQ_r=\calQ$ of $\calQ$ by coherent subsheaves
such that each successive quotient $\calQ_i/\calQ_{i-1}$ is {\it scheme-theoretically} concentrated
on $X$. In particular, we may regard each quotient $\calQ_i/\calQ_{i-1}$ just as a coherent sheaf on $X$.
We set
$$
{\det}_{X}(\calF_1,\calF_2)=\bigotimes\limits_{i=1}^r \det(\calQ_i/\calQ_{i-1}).
$$
It is easy to see that the result does not depend on the choice of the above filtration.

In the general case, there exists another locally free sheaf $\calF_3$ and embeddings
$\beta_1,\beta_2:\calF_1,\calF_2\hookrightarrow\calF_3$ such that

a) Both $\beta_1$ and $\beta_2$ are isomorphisms away from $X$.

b) $\alp=\beta_1|_{S^0}\circ\beta_2|_{S^0}^{-1}$.

\noindent
In this case we set $\det_X(\calF_1,\calF_2,\alp)=\det_X(\calF_1,\calF_3,\beta_1)\ten\det_X(\calF_2,\calF_3,\beta_2)^{-1}$.
It is easy to see that the result is canonically independent of the choice of $\calF_3$ and $\beta_1,\beta_2$ (this definition
is certainly well-known, but we were unable to find a good reference for it).

Let now $G$ be an algebraic group and let $V$ be a finite-dimensional representation of $G$.
In this case for any $G$-bundle $\calF$ (on any variety) we can consider the associated vector bundle
$\calF_V$. Given two $G$-bundles $\calF_1$ and $\calF_2$ on $S$ with an isomorphism between their
restrictions to $S^0$, we set
$$
{\det}_{X,V}(\calF_1,\calF_2)={\det}_X((\calF_1)_V,(\calF_2)_V).
$$

\ssec{}{Representations and bilinear forms}
We now go back to the case when $G$ is reductive and keep all the notations from the previous sections.
Then any finite-dimensional representation $V$ of $G$ as above defines a
symmetric bilinear $W$-invariant form on $\Lam$ in the following way.
By restricting $V$ to $T$ we get a collection of $n=\dim V$ weights $\lam_1^{\vee},...,\lam_n^{\vee}$ of $T$, which
gives an action of $\Lam$ on $\ZZ^n$. We set
$$
Q_V(\lam,\mu)=-\Tr_{\ZZ^n}(\lam\cdot\mu).
$$
It is clear that this form is non-positive definite, integral and $W$-invariant. Also,
if we assume that the determinant of
$V$ is an even character of $G$ (i.e. it is a square of another character) then this form is even
(in particular, that is always the case when $G$ is semi-simple). This form is negative-definite
if $V$ is almost faithful (i.e. if $V$ is a faithful representation of a quotient of $G$ by a finite
central subgroup).

Moreover, it is well-known that for any negative-definite $W$-invariant form $Q$ there exists a positive
integer $c$ and a representation $V$ as above such that $cQ=Q_V$. It is easy to see that
\reft{main} holds for the form $Q$ if and only if it holds for $cQ$. Thus we may assume that
$Q=Q_V$.
\ssec{central-ext}{Description of the central extension}We now want to define the central extension $G_{\aff}$ of the
group $G[t,t^{-1}]\rtimes \GG_m$ (we shall only do it in the case $Q=Q_V$). For any test $\kk$-scheme
$\calW$ we need to define $G_{\aff}(\calW)$. Let us first do it just for
$G[t,t^{-1}]$, i.e. let us define $\tilG(\calW)$.

Consider some element $g\in G[t,t^{-1}](\calW)$; by definition this is the same as a map
$\calW\x\GG_m\to G$. Let $S=\AA^1\x\calW$. Then to $g$ there corresponds a triple
$(\calF_1,\calF_2,\alp_g)$ as above, where $\calF_1$ and $\calF_2$ are both trivial bundles on $\calW\x\AA^1$ and
$\alp_g$ is an isomorphism between their restrictions to $\calW\x\GG_m$ which is "equal to $g$" in the obvious sense.
Then we can consider $\calL_g:=\det_{\calW}(\calF_1,\calF_2)$ which is a line bundle on $\calW$. It is clear that
$\calL_e$ is canonically trivial (here $e$ is the identity element). Moreover, for any
$g_1,g_2\in G[t,t^{-1}](\calW)$ we have the natural isomorphism
\eq{grstr}
\calL_{g_1}\ten\calL_{g_2}\simeq \calL_{g_1g_2}.
\eeq
We set $\tilG(\calW)$ to be the set of pairs $(g,\kap)$ where $\kap$ is a trivialization of $\calL_g$.
The group structure on $\tilG(\calW)$ comes from \refe{grstr}.

It is easy to see that with the above definition of $\tilG$ the action of $\GG_m$ on $G[t,t^{-1}]$
by loop rotations extends to $\tilG$ (this follows from the fact that the line bundle $\calL_g$
considered above does not change when we change $g$ by a loop rotation) and thus we may consider the semi-direct product
$G_{\aff}=\tilG\ltimes \GG_m$.
\ssec{exttors}{Torsors of extension}
Here is another source of examples of $\ZZ$-torsors that will be important in the future. Let $X$ be a curve
and let $x$ be a smooth point of $X$. Set $X^0=X\backslash\{ x\}$. Assume that we are given a line bundle
$\calL^0$ on $X^0$. Then the set of all possible extensions $\calL$ of $\calL^0$ to the whole of $X$ (as a line
bundle) is naturally a $\ZZ$-torsor, which we shall call the torsor of extensions of $\calL^0$ to $X$.

The following remark will become very important later. Assume that $X$ is projective. Then the above torsor of extensions
is canonically trivial. Indeed to any extension $\calL$ we can associate the integer $\deg \calL$ (the degree of $\calL$).
\ssec{}{More $\ZZ$-torsors}
Let now $S$ be a smooth surface and let $X_1,...,X_n$ be a collecton of smooth curves $S$, intersecting
at a point $y\in S$; let us denote their union by $X$. Let now $\calF_1$ and $\calF_2$ be two vector bundles
on $S$ with an isomorphism
away from $X$. Let $X_i^0=X\backslash\{ y\}$.
Then we can consider the relative determinants $\det_{X_i^0}(\calF_1,\calF_2)\in\Pic(X_i^0)$.
Let $\calT_i$  be the $\ZZ$-torsor
of extensions of $\det_{X_i^0}(\calF_1,\calF_2)$ to the whole of $X_i$. The following result is probably well-known, but we were not able to find a reference.
The proof given below is due to V.~Drinfeld.
\lem{drinfeld}The $\ZZ$-torsor
$$
\bigotimes\limits_{i=1}^n\calT_i
$$
is canonically trivial.
\elem
\prf
Clearly, we can replace $S$ by the formal neighbourhood of $y$. Let also
$\calD$ denote the one-dimensional formal disc (formal neighbourhood of
$0$ in $\AA^1$).
Let us choose formal coordinates $u,v$ around $y$. Then for any $a,b\in\kk$
we can consider the map $f_{a,b}:S\to \calD$ sending $(u,v)$ to $au+bv$.
Then for generic pair $(a,b)$ this map induces an isomorphism $X_i\simeq \calD$.
Let us choose such a pair.

As before we may assume that the identification $\calF_1|_{S\backslash X}\simeq \calF_2|_{S\backslash X}$
comes from an embedding $\calF_1\hookrightarrow \calF_2$. Let $\calQ=\calF_2/\calF_1$. This is
a coherent sheaf on $S$, which is set-theoretically supported on $X$. Consider the direct
image $\calR=(f_{a,b})_*\calQ$. Then we may consider $\det \calR$. This is a line bundle on $\calD$;
it is clear that its restriction to $\calD^0$ is canonically isomorphic to
the tensor product of all the $\det_{X_i^0}(\calF_1,\calF_2)$ (where we identify
$X_i^0$ with $\calD^0$ by means of $f_{a,b}$). Thus the $\ZZ$-torsor
$$
\bigotimes\limits_{i=1}^n\calT_i
$$
is canonically isomorphic to the torsor of extensions of $\det\calR|_{\calD^0}$ to $\calD$. But since
$\det\calR$ is a line bundle defined on all of $\calD$, it defines such an extension canonically and
thus our $\ZZ$-torsor is canonically trivial.

The collection of all pairs $(a,b)$ for which the above arguments works is a Zariski open
subset $\calX$ of $\AA^2$; in particular, $\calX$ is an  irreducible algebraic variety.
It is clear that the above trivialization depends regularly on $(a,b)\in\calX$; thus the irreducibility of $\calX$
implies that the trivialization is independent of the choice of the pair $(a,b)$.
\epr
It the sequel, we shall say that a curve $X\subset S$ is {\it good} if it is a union of smooth irreducible components.

\ssec{}{Relative $c_2$}Let again $S$ be a smooth surface and let $X$ be a smooth connected projective curve
inside $S$. Let also $\calF_1,\calF_2$ be two  $G$-bundles on $S$.
Let also $\alp$ be an identification between $\calF_1$ and $\calF_2$ on $S^0=S\backslash X$ (compatible
with the trivialization of the determinants). Then we set
$$
c_2(\calF_1,\calF_2)=\deg({\det}_X(\calF_1,\calF_2)).
$$
More generally, assume that $X$ is a good connected curve. In this case $\det_X(\calF_1,\calF_2)$
is not well-defined as a line bundle on $X$ (it is only well-defined on the smooth part $X^0$ of $X$).
It follows, however, from \refl{drinfeld} that the degree of $\det_X(\calF_1,\calF_2)$ does make sense and we again denote it
by $c_2(\calF_1,\calF_2)$.

Assume now that $G=GL(n)$. In this case we can think about $\calF_1$ and $\calF_2$ as vector bundles of rank $n$ and
we can define the {\em relative Euler characteristic} $\chi(\calF_1,\calF_2)$
in the following way. First, as before, we choose some vector bundle $\calF_3$ on $S$ of rank $n$ containing
both $\calF_1$ and $\calF_2$ as locally free subsheaves. We set $\calQ_i=\calF_3/\calF_i$ for $i=1,2$.
Then both $\calQ_1$ and $\calQ_2$ are set-theoretically supported on $X$ and, therefore, their cohomology is
finite-dimensional. We set
$$
\chi(\calF_1,\calF_2)=\chi(\calQ_1)-\chi(\calQ_2),
$$
where $\chi(\calQ_i)$ denotes the Euler characteristic of $\calQ_i$. It is clear that the result does not depend on the
choice of $\calF_3$.
\prop{chictwo}
Assume that the triple $(\calF_1,\calF_2,\alp)$ is such that $\det(\alp)$ (which is an isomorphism
between $\det\calF_1|_{S^0}$ and $\det\calF_2|_{S^0}$) extends (as an isomoprhism) to the whole of $S$. Then
$$
\chi(\calF_1,\calF_2)=c_2(\calF_1,\calF_2)
$$
\eprop
\prf
Assume first that $X$ is irreducible.
Let choose $\calF_3$ as above; since both $\calQ_1$ and $\calQ_2$ have a filtration with quotients supported on $X$, it follows
that it makes sense to speak about the rank and the degree of $\calQ_i$ on $X$,
which we shall denote by $\rk \calQ_i$ and $\deg\calQ_i$. Then the fact that
$\det(\alp)$ extends to the whole of $X$ implies that $\rk \calQ_1=\rk\calQ_2$. On the other hand, by the definition
we have
$$
c_2(\calF_1,\calF_2)=\deg\calQ_1-\deg\calQ_2.
$$
Thus \refp{chictwo} follows from the Riemann-Roch theorem for the curve $X$.

In the general case (i.e. when $X$ is not necessarily irreducible) we can find a sequence
$(\calF_0=\calF,\calF_1,\cdots,\calF_m=\calF'$ of vector bundles of rank $n$ on $S$
together with the isomorphisms $\alp_i:\calF_i|_{S^0}\widetilde{\to}\calF_{i+1}|_{S^0}$
for all $i=0,\cdots,m-1$ such that

a) The composition of all the $\alp_i$ is equal to $\alp$.

b) Each $\alp_i$ is an isomorphism away from one irreducible component of $X$.

\noindent
In this case the above proof shows that for all $i=0,\cdots,m-1$ we have
$c_2(\calF_i,\calF_{i+1})=\chi(\calF_i,\calF_{i+1})$. Since we have
$$
c_2(\calF,\calF')=\sum\limits_{i=0}^{m-1} c_2(\calF_i,\calF_{i+1})\qquad\text{and}\qquad
\chi(\calF,\calF')=\sum\limits_{i=0}^{m-1} \chi(\calF_i,\calF_{i+1}),
$$
the assertion of \refp{chictwo} follows.
\epr

\refp{chictwo} implies the following result, which in some sense explains the notation
$c_2(\calF,\calF_2)$ (this result will not be used in this paper and we leave the proof to the reader):
\lem{ctwo}
Assume that $S$ is projective. Assume also that $G=SL(n)$ (in other words,
$\calF_1$ and $\calF_2$ are vector bundles, $\det\calF_1$ and $\det\calF_2$ are trivialized and
$\alp$ is compatible with these trivializations). Then
$$
c_2(\calF_1,\calF_2)=c_2(\calF_2)-c_2(\calF_1).
$$

\elem
\ssec{}{The torsor $\calT_{\calF}$}Let now $S$ and $X$ be as above and let $\calF$ be a $G$-bundle on $S^0$. Then to this data we can associate
a 
$\ZZ$-torsor $\calT_{\calF}$, which is uniquely characterized by the following
properties:

1) Any extension $\calF_1$ of $\calF$ to $S$
defines a trivialization of $\kap_{\calF_1}$ of $\calT_{\calF}$. 

2) Let $\calF_1$ and $\calF_2$ be two extensions of $\calF$ to $S$. Then the difference
between $\kap_{\calF_1}$ and $\kap_{\calF_2}$ is equal to $c_2(\calF_1,\calF_2)$.

\ssec{}{The case of disconnected $X$}We would like to refine slightly the above definitions.
Namely, assume $X$ has $r$ connected components. Then
the relative Chern class $c_2(\calF_1,\calF_2)$ naturally takes values in $\ZZ^r$ (since it makes sense to talk about
$c_2(\calF_1,\calF_2)$ in the formal neighbourhood of each connected component of $X$). In particular,
in this case we shall denote by $\calT_{\calF}$ the corresponding $\ZZ^r$-torsor.
\ssec{}{}
Consider the groupoid $\Bun_G(S^0)$. We let $\tilBun_G(S^0)$ denote the "total space of all the $\calT_{\calF}$'s"
(in other words, a point of $\tilBun_G(S^0)$ is given by a $G$-bundle $\calF$ on $S^0$ together with the trivialization
of the torsor $\calT_{\calF}$). By the definition we have a natural map $q:\Bun_G(S)\to \tilBun_G(S^0)$
which is representable.

\th{finite-kapranov}Assume that $X$ can be blown down, i.e. that there exists a
(not necessarily smooth) surface $\Sig$ with a point $y\in \Sig$ and a proper birational
map $f:S\to\Sig$ such that

a) $X$ is equal to the set-theoretic preimage of $y$

b) The restriction of $f$ to $S^0$ is an isomorphism between $S^0$ and $\Sig^0=\Sig\backslash\{y\}$.

\noindent
Assume also that $G$ is semi-simple. Then
\begin{enumerate}

\item
Let $\tcF\in\tilBun_G(S^0)$ and assume that $X$ has $r$ connected components.
Then there exists $A\in\ZZ$ such that for any $\overset{\to}a=(a_1,...,a_r)\in\ZZ^r$ such that
$a_i>A$ for some $i=1,...,r$, we have
$$
\overset{\to}a\cdot \tcF\not\in\im(q)
$$
(here by $\overset{\to}a\cdot \tcF$ we mean the action of $\overset{\to}a$ on $\tcF$ as an element of the
$\ZZ^r$-torsor $\tilBun_G(S^0)$).
\item
$q$ has finite fibers.
\end{enumerate}
\eth

\medskip
\noindent
{\bf Remark.}
When $X$ is smooth this is exactly Theorem 2.2.1 of \cite{Kap-E} (since in this case
our condition is equivalent to the fact that $X$ has negative self-intersection).
\prf
First of all, without loss of generality, we may assume that $X$ is connected; we shall denote by $y\in\Sig$ the image of
$X$ in $\Sig$.

The first assertion of
\reft{finite-kapranov} is equivalent to the following statement. Fix some $\calF\in\Bun_G(S)$. Then we must show
that for all  triples $(\calF,\calF',\alp)$ (where $\calF'\in\Bun_G(S)$ and
$\alp$ is an isomorphism between the restrictions of $\calF$ and $\calF'$ to $S^0$) the set of all possible values
of $c_2(\calF,\calF')$ is bounded above.

Similarly, the second assertion of \reft{finite-kapranov} says that for any $\calF\in\Bun_G(S)$ and $l\in \ZZ$,
the set of isomorphism classes of
triples $(\calF,\calF',\alp)$ as above such that
\eq{ctwozero}
c_2(\calF,\calF')=l
\eeq
is finite.

Since $G$ is semi-simple we can choose an embedding $G\hookrightarrow SL(n)$; without loss of generality we may assume that
the quadratic form $Q$ comes from this embedding. Then the set of all possible $(\calF',\alp)$ for $G$ embeds into the similar
set for $SL(n)$. Hence we may assume that $G=SL(n)$.
In this case $\calF$ and $\calF'$ above should
be thought of as vector bundles of rank $n$ with trivial determinant.

Let us first concentrate on \reft{finite-kapranov}(2). Then
$$
\chi(\calF,\calF')=c_2(\calF,\calF')=0.
$$

The proof of the required finiteness is based on the following result:
\lem{E}There exists a locally free subsheaf $\calE$ of $\calF$, which is equal to $\calF$ on $S^0$ and such that
for any $\calF'$ as above the resulting identification $\calE|_{S^0}\simeq \calF'|_{S^0}$ extends to an embedding
$\calE\hookrightarrow\calF'$.
\elem
\prf
Let $j$ denote the (open) embedding of $S^0$ into $\Sig$. Then we can consider the sheaf $\calG=j_*(\calF|_{S^0})$
on $\Sig$. Since the complement of $S^0$ in $\Sig$ consists of one point, it follows that this sheaf is coherent.
Also, for any $\calF'$ as above we have the natural embedding (of coherent sheaves) $f_*\calF'\hookrightarrow \calG$
which is an isomorphism away from $y$. Hence there exists some $N>0$ such that $f_*\calF'$ contains
$\grm_{\Sig,y}^N\calG$, where $\grm_{\Sig,y}$ denotes the maximal ideal of the point $y$ in $\Sig$.

We claim now that such the number $N$ as above can be chosen uniformly for all $\calF'$ satisfying \refe{ctwozero}.
To prove this let us define
$$
\eta(\calF')=\length(\coker(f_*\calF'\to \calG))+\length(R^1f_*\calF').
$$
Then it is clear that $l=\chi(\calF,\calF')=\eta(\calF)-\eta(\calF')$. Thus
$$
\length(\coker(f_*\calF'\to \calG))=\eta(\calF)-\length(R^1 f_*(\calF'))-l.
$$
Hence $\length(\coker(f_*\calF'\to \calG))\leq \eta(\calF)-l$. Hence (by Nakayama lemma) if we choose a number
$N$ such that $\length(\calG/\grm_{\Sig,y}^N\calG)\geq \eta(\calF)-l$,
then the sheaf
$f_*\calF'$ must contain $\grm_{\Sig,y}^{N}\calG$.

Let now $\calH=\grm_{\Sig,y}^N\calG$. This a coherent sheaf on $\Sig$. The embedding $\calH\hookrightarrow f_*\calF'$
gives rise to a morphism $f^*\calH\to \calF'$, which is an isomorphism on $S^0$. Let $\calE'$ denote the quotient
of $f^*\calH$ by torsion. This is a torsion-free coherent sheaf on $S$. Let $\calE$ denote its saturation
(i.e. minimal locally free sheaf, containing $\calE'$). Then the morphism $f^*\calH\to \calF'$ gives rise
to an embedding $\calE\hookrightarrow \calF'$ of locally free sheaves, which is an isomorphism on $S^0$.
\epr
Let us now explain why \refl{E} implies \reft{finite-kapranov}(1). Without loss of generality, we may assume that
$\calE=\calF(N\cdot X)$. Applying the same argument to
the dual of $\calF$ we see that we may also assume that any $\calF'$ as above is contained in $\calF(-N\cdot X)$.
Let $\calQ=\calF(-N\cdot X)/\calF(N\cdot X)$. This is a coherent sheaf on $S$, which is set-theoretically supported
on $X$. Any $\calF'$ as above is uniquely determined by its image $\calF'_{\calQ}=\calF'/\calF(-N\cdot X)$ in $\calQ$.
Moreover, the Euler characteristic $\chi(\calF'_{\calQ})$ is independent of $\calF'$ (it is equal to
$\chi(\calF,\calF')+\chi(\calF(N\cdot X),\calF)=l+\chi(\calF(N\cdot X),\calF)$). 

Hence it is enough for us to show that there exists a scheme of finite type over $\kk$ whose $\kk$ points are
subsheaves of $\calQ$ with given Euler characteristic (indeed, since $\kk$ is finite, the set of $\kk$-points of a scheme of finite
type is finite). Without loss of generality we can assume that $S$ is projective. We shall choose an embedding of $S$ into some projective space; in other
words we are going to choose some very ample line bundle $\calL$ on $S$. Thus $Q$ is a coherent sheaf on a projective
scheme $S$ and hence according to \cite{Groth} there exists a scheme $Quot(\calQ)$ whose $\kk$-points are subsheaves $\calG$ of $\calQ$. Moreover, the subscheme
of $Quot(\calQ)$ corresponding to looking at $\calG$ as above with fixed Hilbert polynomial is of finite type over $\kk$. For every irreducible component
$X_i$ of $X$ we can define the generic rank $\rk_i(\calG)$ of $\calG$ as follows: first, we can fined a filtration
$0\subset \calG_1\subset \cdots\subset \calG_s=\calG$ of $\calG$ such that every quotient
$\calG_j/\calG_{j-1}$ is scheme-theoretically concentrated on $X$. Then 
we define $\rk_i(\calG)=\sum_j \rk_i(\calG_j/\calG_{j-1})$ where $rk_i(\calG_j/\calG_{j-1})$ is the rank of the restriction
of $\calG_j/\calG_{j-1}$ to the generic point of $X_i$. It is clear that $0\leq \rk_i(\calG)\leq \rk_i(\calQ)$, hence there are finitely many
possibilities for each of the numbers $\rk_i(\calG)$. Thus to finish the 
proof it is enough to show the following
\lem{}The Hilbert polynomial of $\calF$ is uniquely determined by $\chi(\calG)$ and by the numbers $\rk_i(\calG)$.
\elem
\prf
For each $X_i$ as above let us set $\calL_i=\calL|_{X_i}$; let $d_i=\deg _{X_i}\calL_i$. Also, let us choose a filtration $\calG_j$ of $\calG$ as above such that 
every the quotient $\calG_j/\calG_{j-1}$ is concentrated scheme-theoretically on some irreducible component $X_{i_j}$ of $X$; we shall
denote by $\calS_j$ the corresponding coherent sheaf on $X_{i_j}$. By the definition, the Hilbert polynomial of $\calG$ is determined
by the numbers $\dim H^0(\calG\ten\calL^{\ten n})$ for $n$ sufficiently large. However, if $n$ is large then $H^p(\calG\ten\calL^{\ten n})=0$ for $p>0$ and thus
we have 
$$
\begin{aligned}
\dim H^0(\calG\ten\calL^{\ten n})=\chi(\calG\ten\calL^{\ten n})=\sum\limits_{j}\chi(\calS_j\ten \calL_{i_j}^{\ten n})=\\
\sum\limits_{j}\chi(\calS_j)+\rk(\calS_j)\cdot n\cdot d_{i_j}=\chi(\calG)+\sum\limits_{i}\rk_i(\calG)\cdot n\cdot d_i,
\end{aligned}
$$
which finishes the proof.
\epr

\smallskip
\noindent
{\bf Remark}. The same proof works for arbitrary reductive $G$ if we fix the $G/[G,G]$-bundle obtained from $\calF'$
by push-forward with respect to the natural map $G\to G/[G,G]$.

\smallskip
\noindent
Let us now prove \reft{finite-kapranov}(1). We are going to use the notations introduced in the proof of \refl{E}.
We need to show that $\chi(\calF,\calF')$ is bounded below. However, we have
$$
\chi(\calF,\calF')=\eta(\calF)-\eta(\calF').
$$
By definition, $\eta(\calF')\geq 0$. Hence $\chi(\calF,\calF')\leq \eta(\calF)$.
\epr


\sec{central}{Turning on the central extension}
\ssec{}{}We now want to generalize \refp{bijectionw} to the case where
the group $G_{\aff}'$ is replaced by its central extension using the results of the previous Section.
For this we have to first interpret the $\ZZ$-torsor $G_{\aff}(\calO)\backslash \tilpi^{-1}(k)/G_{\aff}(\calO)$
over $G'_{\aff}(\calO)\backslash \pi^{-1}(k)/G'_{\aff}(\calO)$ in geometric terms. To do this
we are going to apply the constructions of the previous Section to $S=\tilS_k$ where the role of $X$ will
be played by  the exceptional fiber of the
morphism $\tilS_k\to S_k$ which we shall denote by $E$ (we want to reserve the notation $X$ for something else).
Note that in this case $S^0=S_k$; we shall write $\tilBun_G(S^0_k)$ instead
of writing $\tilBun_G(S^0)$.
\prop{centext}
The groupoid
$G_{\aff}(\calO)\backslash \tilpi^{-1}(k)/G_{\aff}(\calO)$ is canonically equivalent
to the groupoid $\tilBun_G(S^0_k)$.
\eprop
\prf
Consider the quotient $\tilG(\calK)/\calO^*$ (where we take the quotient by the central $\calO^*\subset\calK^*$).
This is a central extension of $G(\calK[t,t^{-1}])$ by $\ZZ$. In particular, it defines a $\ZZ$-torsor
over every $g(t,s)\in G(\calK)$ which we shall denote by $\calT_g$. By the definition, this torsor  can be described as follows.

Consider the
plane $\AA^2$ with coordinates $t$ and $s$. Let $C\subset \AA^2$
be given by the equation $s=0$ and let $D\subset \AA^2$ be given by the equation $t=0$. Let $y\in\AA^2$ be the point
$t=0,s=0$. Set $C^0$ (resp. $D^0$) to be the complement of $y$ in $C$ (resp. $D$).
Let now $\calM_1,\calM_2$ be two trivial bundles on $\AA^2$. Let $\xi$ be the isomorphism
between their restrictions on $\GG_m\x\GG_m$ given by $g$ ($\xi$ is the natural isomorphism between two trivial bundles on
$\GG_m\x\GG_m$ multiplied by $g(t,s)$). Consider $\det_{D^0}(\calM_1,\calM_2)$; this is a line bundle on $D^0$ (which is
isomorphic to $\GG_m$ with coordinate
$s$). Then by the definition, the $\ZZ$-torsor $\calT_g$ is the torsor of extensions of
$\det_{D^0}(\calM_1,\calM_2)$ from $D^0$ to $D$.
On the other hand, by \refl{drinfeld} this torsor is canonically isomorphic to the inverse of the torsor of extensions
of $\det_{C^0}(\calM_1,\calM_2)$ to $C$.

On the other hand, given $g(t,s)$ as above we can construct a $G$-bundle $\calF$ on $S^0_k$ together with trivializations
on $U_k$ and $V_k$ as in the proof of \refp{bijectionw}. What we need to do, is to construct a trivialization
of $\calT_g$ starting from a trivialization of $\calT_{\calF}$ (recall that we consider $S_k^0$ as an open subset
of $\tilS_k$ and the torsor $\calT_{\calF}$ is constructed via this embedding).
In other words, we need to start with an extension
$\calF'$ of $\calF$ to $\tilS_k$ and construct a trivialization of $\calT_g$ (the fact that this construction gives
rise to an isomorphism $\calT_{\calF}\simeq \calT_g$ will eventually be obvious).

Let $X\subset S_k$ be given by the equation $x=0$ in $S_k$; let $X^0=X\backslash\{ 0,0,0\}=S_k^0\cap X$.
Let $\alp$ denote the trivialization of $\calF$ on $U_k=\Spec\kk[x,x^{-1},y]$. Let us also identify
$X^0$ with $C^0$ by setting $y=t^{-1}$. Then it is easy to see that $\det_{C^0}(\calM_1,\calM_2)=\det_{X^0}(\alp)$.
Thus the torsor of extensions of $\det_{C^0}(\calM_1,\calM_2)$ to $0$ is canonically isomorphic to the torsor of
extensions of $\det_{X^0}(\alp)$ to $\infty$. The latter torsor is inverse to the torsor of extensions
of $\det_{X^0}(\alp)$ to $0$. In other words we get the canonical isomorphism between $\calT_g$ and the torsor of extensions
of $\det_{X^0}(\alp)$ to $0$.

We claim now, that an extension $\calF'$ of $\calF$ to $\tilS_k$
as above provides a canonical extension of $\det_{X^0}(\alp)$ to
the whole of $X$. Indeed, let us identify $X$ with its proper
preimage in $\tilS_k$ and let also $E$ denote the exceptional
divisor in $S_k$. Then $E$ and $X$ intersect at a unique point $p$
in $\tilS_k$ and $Y=E\cup X$ is a good curve in $\tilS_k$; we set
$E^0=E\backslash\{p\}$. Let $\calF^0$ denote the trivial bundle on
$\tilS_k$; we may regard $\alp$ as an isomorphism between
$\calF'|_{S_k^0}$ and $\calF^0|_{S_k^0}$. Then by \refl{drinfeld}
the torsor of extensions of $\det_{X^0}(\alp)$ to $X$ is inverse
to the torsor of extensions of $\det_{E^0}(\calF',\calF^0,\alp)$
to $E$. However, since $E$ is proper, by \refss{exttors} the latter torsor acquires a
canonical trivialization and hence the same is true for the
former. \epr

\ssec{}{Geometric interpretation of the convolution diagram}
We now want to generalize the results of \refss{convolution} in order to give a geometric interpretation of the convolution
diagram in the presence of a central extension. Namely, let us choose as before two positive integers
$k$ and $l$. Recall that in this case we have the partial resolution $\grp_{k,l}:S_{k,l}\to S_{k+l}$.
Also we have the full resolution $\grp_{1,\cdots,1}:\tilS_{k+l}=S_{1,\cdots,1}\to S_{k+l}$ which factorizes
through a map $\grr:\tilS_{k+l}\to S_{k,l}$ composed with $\grp_{k,l}$.

The map $\grr$ is an isomorphism over $S^0_{k,l}$. Note that the complement of $S^0_{k,l}$ in $\tilS_{k+l}$
is disconnected (it has two connected components, which correspond to the two singular points of $S_{k,l}$).
Hence the groupoid $\tilBun_G(S^0_{k,l})$ (here we view $S^0_{k,l}$ as a complement to a divisor in the smooth
surface $\tilS_{k+l}$)
is $\ZZ\x\ZZ$-torsor over $\Bun_G(S^0_{k,l})$. Since $S_{k,l}$ contains $S_{k+l}$ as an open subset, it follows
that the restriction map $\Bun_G(S^0_{k,l})\to \Bun_G(S^0_{k+l})$ extends to a map
$$
\tilBun_G(S^0_{k,l})\to \tilBun_G(S^0_{k+l}).
$$
This map is $\ZZ\x\ZZ$-equivariant, where $\ZZ\x\ZZ$ acts naturally on the left hand side and
on the right hand side it acts by means of the homomorphism $\ZZ\x\ZZ\to\ZZ$
sending $(a,b)$ to $a+b$.

\prop{conv-central}The groupoid
\eq{long}
G_{\aff}(\calO)\backslash\tilpi^{-1}(k)\underset{G_{\aff}(\calO)}\x \tilpi^{-1}(l)/G_{\aff}(\calO)
\eeq
can be canonically identified with $\tilBun_G(S^0_{k,l})$. Under this identification the natural
$\ZZ\x\ZZ$-action on \refe{long} corresponds to the $\ZZ\x\ZZ$-action on $\tilBun_G(S^0_{k,l})$.

Moreover, under this identification the multiplication map
$$
G_{\aff}(\calO)\backslash\tilpi^{-1}(k)\underset{G_{\aff}(\calO)}\x \tilpi^{-1}(l)/G_{\aff}(\calO)\to
G_{\aff}(\calO)\backslash \tilpi^{-1}(s^{k+l})/G_{\aff}(\calO)
$$
corresponds to the map $\tilBun_G(S^0_{k,l})\to\tilBun_G(S^0_{k+l})$
discussed above.
\eprop
The proof of \refp{conv-central} is essentially the same as the proof of \refp{centext} and we leave it to the reader.
\ssec{}{Proof of \reft{main}(1)}
We now claim that \refp{conv-central} together with \reft{finite-kapranov} imply \reft{main}(1).
First of all, this is clear if $G$ is semi-simple: in this case condition 1') of \refs{general} follows immediately
from \reft{finite-kapranov}(1) and condition 2) follows from \reft{finite-kapranov}(2). Also, when $G$ is a torus
\reft{main}(1) follows from \refs{ex}. Hence it is easy to see that \reft{main}(1) holds true for any reductive group
$G'$ which is isomorphic to a product of a semi-simple group and a torus.
Also assume that we are given an isogeny $G\to G'$ of split reductive groups. Then it is easy to see that \reft{main}(1)
for $G'$ implies \reft{main}(1) for $G$. By taking $G'$ to be the product of the adjoint group of $G$ and
of $G/[G,G]$, we see that \reft{main}(1) is true for any $G$.
\sec{principal}{Action on the principal series and proof of
\reft{main}(2)}
\ssec{}{Principal series}Consider first the subgroup $\tilDel'$ of $\tilGam$ which is equal  to
$\calO^*\x G(\calK[t])\rtimes \calO^*$. We shall denote by $\Del'$ its image in $\Gam$.

Now we define subgroups $\Del\subset \Del'$ and $\tilDel\subset \tilDel'$. First of all note that $\tilDel'$
maps naturally
to $\calO^*\x G(\calK)\x \calO^*$ (by setting $t=0$). Let $B$ be a Borel subgroup of $G$. It has the decomposition
$B=TU$, where $T$ is a maximal torus in $G$ and $U$ is the unipotent radical of $B$.
We set $\tilDel$  to be the preimage of $\calO^*\x T(\calO)U(\calK)\x\calO^*$ in $\tilDel'$.
We denote by $\Del$ the image of $\tilDel$ in $\Gam$.

We now set
$$
\Ome=\Gam/\Del,\  \tilOme=\tilGam/\tilDel,\  \Ome'=\Gam/\Del',\  \tilOme'=\tilGam/\tilDel'.
$$
Clearly, $\Gam$ acts on $\Ome$ and $\Ome'$ and $\tilGam$
acts on $\tilOme$ and $\tilOme'$. In addition we have the natural map
$\varpi:\tilOme\to\tilOme'$, which is $\tilGam$-equivariant.
\th{goodaction}
The action of $\tilGam^+$ on $\tilOme$ is good with respect to $\tilOme'$.
\eth
\ssec{}{The structure of $\Ome$}Before we go to the proof of \reft{goodaction} we would like to discuss how
$\tilOme$ and $\tilOme'$ look like.
First of all, let us introduce another pair of groups $(\Del'',\tilDel'')$ such that
$\Del\subset\Del''\subset\Del'$ (resp. $\tilDel\subset\tilDel''\subset\tilDel'$). These are defined exactly
as $\Del$ and $\tilDel$ except that the group $T(\calO)U(\calK)$, used in the definition, has to be replaced by
$B(\calK)$. It is easy to see that $\Del$ is normal in $\Del''$ (resp. $\tilDel$ is normal in $\tilDel''$)
and we have the natural isomorphism $\Del''/\Del=\tilDel''/\tilDel=T(\calK)/T(\calO)\simeq \Lam$.
Thus $\Lam$ acts on $\Ome$ on the right. In addition the group
$\ZZ$ acts on $\Ome$ (this action comes from the center of $\Gam$). Altogether, we get an action of $\ZZ\x\Lam=\Lam_{\aff}'$
on $\Ome$ and this action commutes with the action of $\Gam$. In particular, $\Lam_{\aff}'$ acts on $\Gam_0\backslash\Ome$.
Similarly, the lattice $\Lam_{\aff}=\ZZ\x\Lam\x\ZZ$ acts on $\tilOme$ and this action commutes with the action of
$\tilGam$. In particular, $\Lam_{\aff}$ acts on $\tilGam_0\backslash\tilOme$.
The following lemma follows easily from Proposition 1.4.5 in \cite{Kap-H}.
\lem{strome}
\begin{enumerate}
\item
The action of $\Lam_{\aff}'$ on $\Gam_0\backslash\Ome$ is simply transitive.
Similarly, the action of $\Lam_{\aff}$ on $\tilGam_0\backslash\tilOme$ is simply transitive.
\item
 The action of  $\ZZ\x \ZZ$ on $\tilGam_0\backslash \tilOme'$ is simply transitive.
\item
Let us identify $\Gam_0\backslash\tilOme$ with $\Lam_{\aff}=\ZZ\x \Lam\x\ZZ$ by acting on the unit element.
Similarly, let us identify $\tilGam_0\backslash\tilOme'$ with $\ZZ\x\ZZ$. Then the natural
map $\tilGam_0\backslash\tilOme\to\tilGam_0\backslash\tilOme'$ induced by $\varpi:\tilOme\to\tilOme'$ corresponds
to the projection onto the first and the third factor.
\end{enumerate}
\elem
\refl{strome} implies that $\calF^{fin}(\tilOme)\simeq \CC[\Lam_{\aff}]$. Similarly,
$\calF_{\tilOme'}(\tilOme)$ is the $v$-adic completion of $\CC[\Lam_{\aff}]$ (here as before the
action of $v$ on the right hand side comes from the ``central'' copy of $\ZZ$ in $\Lam_{\aff}$)), which consists of all
functions $f(k,\lam,n)$ such that

1) $f(k,\lam,n)=0$ for $k\ll 0$

2) For almost all $n$ we have $f(k,\lam,n)=0$ for all $\lam,k$.

3) For given $k$ the function $f(k,\lam,n)$ is not equal to $0$ only for finitely
many pairs $(\lam,n)$.
\ssec{}{$\tilOme'$ is almost good}
Let us now move to the proof of \reft{goodaction}. By the definition we first have to prove that the action
of $\tilGam$ on $\tilOme'$ is almost good.
The proof is similar to that of \reft{main}. Namely, we are going to give a geometric interpretation of the
convolution diagram for the action of $\tilGam$ on $\tilOme'$. More precisely, for any $l\in\ZZ$, let
$\tilOme'_l=\tilpi^{-1}(l)/\tilDel'$ (resp. $\Ome'_l=\pi^{-1}(l)/\Del'$). Then for any $k>0$ we want to give a
a geometric interpretation of the morphism (of groupoids)
$$
\tilGam_0\backslash\tilOme'_l\leftarrow
\tilGam_0\backslash\tilpi^{-1}(k)\underset{\tilGam_0}\x\tilOme'_l\to \Gam_0\backslash\tilOme'_{k+l}.
$$

First let us do for $\Gam$ and $\Ome'$. In other words,
we want to give a geometric interpretation of the morphism (of groupoids)
$$
\Gam_0\backslash\Ome'_l\leftarrow
\Gam_0\backslash\pi^{-1}(k)\underset{\Gam_0}\x\Ome'_l\to \Gam_0\backslash\Ome'_{k+l}.
$$
First, we need some preparatory material.
\ssec{}{Blow ups}Let $\Sig$ be a smooth surface and let $Y$ be a smooth curve in $\Sig$. Let $x$ be some point in $Y$.
Set $\Sig^0=\Sig\backslash\{ x\}$.
Let $\calO_{\Sig}$ (resp. $\calO_Y$) denote the structure sheaf of $S$ (resp. of $Y$) and let
$\grm_{Y,x}$ be the ideal sheaf of the point $x$ in $\calO_Y$. Abusing the notation, we shall
regard  $\calO_Y$ as a sheaf on $S$ with support on $Y$.
For each $k>0$ let $\calJ_k$ be the preimage of $\grm_{Y,y}^k$ under the natural map $\calO_{\Sig}\to\calO_Y$.
We shall denote by $\Bl_{x,k}(\Sig)$ the blow up of $\Sig$ at the ideal $\calJ_k$. This is a new surface which
has unique singular point $y$ of type $A_k$
\footnote{in fact, if $k=1$, then this point is smooth; in this case $\Bl_{x,k}(\Sig)$ is the usual blow up
of $S$ at the point $x$ and the point $y$ corresponds to the tangent space to $Y$ at $x$}
and it is endowed with a proper birational map $f:\Bl_{x,k}(\Sig)\to \Sig$,
which is an isomorphism away from $x$ and whose fiber over $x$ is isomorphic to $\PP^1$. We shall denote
by $\Bl_{x,k}(\Sig)^0$ the complement to the point $y$ in
$\Bl_{x,k}(\Sig)$. Also we shall identify the proper preimage of $Y$ in $\Bl_{x,k}(\Sig)$ with $Y$. Let
$\Sig'=\Bl_{x,k}(\Sig)\backslash Y$.
Thus we have the natural morphism of groupoids $\grq:\Bun_G(\Bl_{x,k}(\Sig)^0)\to \Bun_G(\Sig)$ and
$\grp:\Bun_G(\Bl_{x,k}(\Sig)^0)\to \Bun_G(\Sig')$
(the first morphism is obtained by restricting a $G$-bundle to $\Sig^0$ and then taking the unique extension to $\Sig$;
the second
morphism is obtained just by restriction to $\Sig'$).

The surface $\Bl_{x,k}(\Sig)$ admits canonical resolution of singularities $\tBl_{x,k}(\Sig)$. This is a smooth surface,
containing $\Bl_{x,k}(\Sig)^0$ as an open subset; moreover, the complement to
$\Bl_{x,k}(\Sig)^0$ in $\tBl_{x,k}(\Sig)$ is a good connected projective curve.
Thus it makes sense to consider $\tilBun_G(\Bl_{x,k}(\Sig)^0)$.
We have canonical morphism $\tgrq:\tilBun_G(\Bl_{x,k}(\Sig)^0)\to \Bun_G(\Sig)\x \ZZ$  (which modulo $\ZZ$
give rise to $\grq$) constructed as follows. Let us pick up an element in
$\tilBun_G(\Bl_{x,k}(\Sig))$. We may assume that it comes from a $G$-bundle $\tcF$ on $\tBl_{x,k}(\Sig)$.
To describe $\tgrq(\tcF)$ we only need to describe its projection to $\ZZ$.
Consider $\grq(\tcF|_{\Bl_{x,k}(\Sig)^0})$. Let us denote by $\calF'$ its lift to $\tBl_{x,k}(\Sig)$.
Then the desired integer is $c_2(\calF,\calF')$.

We shall also denote by $\tgrp:\tilBun_G(\Bl_{x,k}(\Sig)^0)\to \Bun_G(\Sig')$ the composition
of the projection $\tilBun_G(\Bl_{x,k}(\Sig)^0)\to \Bun_G(\Bl_{x,k}(\Sig)^0)$ with $\grp$.

Here is an example of such a situation: let $X$ be another smooth curve and let $\calL$ be a line bundle over $X$. Let also $x$ be
some point of $X$.
We want to consider the above construction in the case when $S=\calL$ and $Y=\calL_x$ -- the fiber of $\calL$ over
$x$. In this case it is easy to see that $\Bl_{x,k}(\Sig)\backslash Y$ is naturally isomorphic to $\calL(-kx)$.

Let us now consider the case $X=\AA^1$. In this case the bundle $\calL$ is automatically trivial
and thus we can identify $\Sig$ with $\AA^2$ with coordinates $(p,s)$, where the point $x$ corresponds to $(0,0)$ and the
curve $Y$ is given by the equation $s=0$.
In this case the surface $\Bl_{x,k}(\Sig)$ can be described as follows.
Consider two surfaces $\Sig_1$ and $\Sig_2$, where $\Sig_1=\Spec\kk[p,q,s]/pq-s^k$ and
$\Sig_2=\Spec\kk[u,s]$ (here $p,q$ and $u$ are some additional variables). Let us identify the open subset
of $\Sig_1$, given by the equation $q\neq 0$, with the open subset of $\Sig_2$, given by the equation
$u\neq 0$ by setting $u=q^{-1}$. Then $\Bl_{x,k}(\Sig)$ is obtained by gluing $\Sig_1$ and $\Sig_2$ along this common
open subset. The corresponding map $\Bl_{x,k}(\Sig)\to\Sig$ is obvious: on $\Sig_1$ it is given
by $(p,q,s)\mapsto (p,s)$ and on $\Sig_2$ it is given by $(u,s)\mapsto (us^k,s)$. The proper preimage of $Y$ lies inside
$\Sig_1$ and is given there by the equation $q=0$. Thus in this case $\Sig'=\Sig_2$ and it can be naturally
identified with $\AA^2$ with coordinates $(u,s)$.

In what follows we shall only work with $\Sig,\Sig'$ in this case. We shall also set
$\Sig_1^0=\Sig_1\cap \Bl_{x,k}(\Sig)^0$. We have natural isomorphisms $\Sig_1\simeq S_k$ and $\Sig_1^0\simeq S_k^0$.

Now we claim the following
\prop{conv-act}The diagram of groupoids
$$
\Gam_0\backslash\Ome'_l\leftarrow
\Gam_0\backslash\pi^{-1}(k)\underset{\Gam_0}\x\Ome'_l\to \Gam_0\backslash\Ome'_{k+l}.
$$
is equivalent to the diagram
$$
\Bun_G(\Sig')\overset{\grp}\leftarrow\Bun_G(\Bl_{x,k}(\Sig)^0)\overset{\grq}\to\Bun_G(\Sig).
$$
Similarly, the diagram
$$
\tilGam_0\backslash\tilOme'_l\leftarrow
\tilpi^{-1}(k)\underset{\Gam_0}\x\tilOme'_l\to \Gam_0\backslash\tilOme'_{k+l}
$$
is equivalent to
$$
\Bun_G(\Sig')\x \ZZ\overset{\tgrp\x\id}\leftarrow\tilBun_G(\Bl_{x,k}(\Sig)^0)\x\ZZ\overset{\tgrq+\id}\to
\Bun_G(\Sig)\x\ZZ
$$
(here by $\tgrq+\id$ we mean the map which sends a pair $(\tcF,a)\in \tilBun_G(\Bl_{x,k}(\Sig)^0)\x\ZZ$,
such that $\tgrq(\tcF)=(\calG,b)$, to $(\calG,a+b)$).
\eprop
\prf
The proof is very similar to the proof of \refp{conv-no-cent} and \refp{conv-central}.
First, for any $n\in\ZZ$ let us consider the quotient
$$
\Gam_0\backslash\Ome'_n=G[ts^n,t^{-1}s^{-n},s]\backslash G[t,t^{-1},s,s^{-1}]/G[t,s,s^{-1}].
$$
We claim that
(as a groupoid) it can naturally be identified with $\Bun_G(\AA^2)$ with coordinates $(u,s)$ where $u=ts^n$.
The proof is exactly the same as the proof of the bijection $(1)\Leftrightarrow (2)$ from \refp{bijectionw}.

Next, we claim that the $\ZZ$-torsor $\tilGam_0\backslash\tilOme'_n\to\Gam\backslash\Ome'_n$ is canonically trivialized.
In other words, we get an equivalence
$$
\tilGam_0\backslash\tilOme'_n\simeq \Bun_G(\Sig)\x \ZZ.
$$
This would have been obvious if we were dealing with sets, rather than with groupoids, since
the set of isomorphism classes of points of $\Gam\backslash\Ome'_n$ consists of one element.
In order to construct the required trivialization on the level of groupoids we must prove the following
result:
\lem{}
The group of automorphisms of the unique point of $\Gam\backslash\Ome'_n$ acts trivially on $\ZZ$.
\elem
\prf
Note, that by definition this group is equal to $\Del'=G(\calK[t])\ltimes \calO^*$.
It is easy to see that it is enough to prove this separately in the case when $G$ is a torus
and when $G$ is semi-simple. If $G$ is a torus then our statement follows from the calculations
of \refs{ex}. When $G$ is semi-simple, then $\Del'$ is actually the group of $\kk$-points of a connected
group ind-scheme $\calG_{\Del'}$ over $k$ and it is clear that the above homomorphism from $\Del'\to\ZZ$
comes from an algebraic homomorphism $\calG_{\Del'}\to\ZZ$, which has to be trivial since $\calG_{\Del'}$
is connected.
\epr

Let us now recall that $\Bl_{x,k}(\Sig)^0$ is glued from $\Sig_1^0\simeq S_k^0$
and $\Sig_2\simeq \AA^2$. Thus setting $p=ts^{k+n},q=ts^{-n}$ and arguing as in the proof
of \refp{conv-no-cent} we get the the equivalence
$$
\Gam_0\backslash\pi^{-1}(k)\underset{\Gam_0}\x\Ome'_l\simeq \Bun_G(\Bl_{x,k}(\Sig)^0).
$$
The rest of the proof is essentially a word-by-word repetition of the proof of \refp{conv-no-cent} and
\refp{conv-central} and it is left to the reader.
\epr
\ssec{delstr}{Proof of \reft{goodaction}}We can now prove \reft{goodaction}.
Indeed, the fact that the action of $\tilGam$ on $\tilOme'$ is good
follows from \refp{conv-act} together with \reft{finite-kapranov} (as in the proof of \reft{main}(1)).
What is left to show is property (i) from \refss{goodactgen}. For this we need some additional geometric construction.

Recall that $\Sig$ is a (trivial) line bundle over $X\simeq \AA^1$. In particular, $X$ is embedded into $\Sig$.
In coordinates $(p,s)$ it is given by the equation $s=0$.

Let now $\Bun_{G,\Del}(\Sig)$ denote the groupoid of $G$-bundles $\calF$ on $\Sig$ endowed with the following additional
data:

1) A reduction of $\calF|_{X^0=X\backslash\{ x\}}$ to $B$. We shall denote the corresponding $B$-bundle on $X^0$
by $\calF_B^0$. We shall also denote by $\calF_T^0$ the corresponding $T$-bundle.

2) An extension $\calF_T'$ of $\calF_T^0$ to $X$.

\smallskip
\noindent
Then arguing as in the previous subsection, for any $n\in\ZZ$ we can  construct an equivalence of groupoids
$$
\Gam_0\backslash\Ome_n\simeq \Bun_{G,\Del}(\Sig),\qquad
\tilGam_0\backslash\tilOme_n\simeq \Bun_{G,\Del}(\Sig)\x\ZZ.
$$
Recall, that the set of isomorphism classes of the groupoid $\Gam_0\backslash\Ome_n$ is $\Lam$.
Geometrically the corresponding map $\Bun_{G,\Del}(\Sig)\to \Lam$ is constructed as follows.
Since the quotient $G/B$ is proper, any $B$-structure on $X^0$ extends uniquely to $X$. Let us denote
by $\calF_B$ the corresponding $B$-bundle on $X$ and by $\calF_T$ the induced $T$-bundle
In particular, this shows that $\Bun_{G,\Del}(\Sig)$ classifies triples $(\calF,\calF_B,\calF_T')$ where
$\calF_B$ is a reduction of $\calF$ to $B$ on $X$ and $\calF_T'$ is a modification of $\calF_T$ at
$x$, i.e. $\calF_T'$ is a $T$-bundle on $X$ endowed with an isomorphism $\calF_T'|_{X^0}\simeq \calF_T|_{X^0}$.
It is well known that the set of isomorphism classed of
pairs of $T$-bundles $(\calF_T,\calF_T')$ on (any smooth curve) $X$ together with an isomorphism on $X^0$
is in one-to-one correspondence with elements of $\Lam$.
Given some $(\calF,\calF_B,\calF_T')\in \Bun_{G,\Del}(\Sig)$ we shall denote by
$d(\calF,\calF_B,\calF_T')$ the corresponding element of $\Lam$ and we shall call it
{\it the defect} of $(\calF,\calF_B,\calF_T')$.

Similarly, we can define $\Bun_{G,\Del}(\Bl_{x,k}^0)$ (in that case we should work with the proper preimage of $X$
in $\Bl_{x,k}(\Sig)$ which we shall identify with $X$) and
$$
\tilBun_{G,\Del}(\Bl_{x,k}^0)=\Bun_{G,\Del}(\Bl_{x,k}^0)\underset{\Bun_{G}(\Bl_{x,k}^0)}\x
\tilBun_{G}(\Bl_{x,k}^0)
$$
together with natural equivalences
$$
\Gam_0\backslash\pi^{-1}(k)\underset{\Gam_0}\x\Ome_l\simeq \Bun_{G,\Del}(\Bl_{x,k}(\Sig)^0),\qquad
\tilGam_0\backslash\tilpi^{-1}(k)\underset{\tilGam_0}\x\tilOme_l\simeq \tilBun_{G,\Del}(\Bl_{x,k}(\Sig)^0).
$$
The corresponding maps
$$
\Bun_{G,\Del}(\Sig')\overset{\grp_{\Del}}\leftarrow\Bun_{G,\Del}(\Bl_{x,k}(\Sig)^0)
\overset{\grq_{\Del}}\to\Bun_{G,\Del}(\Sig)
$$
and
$$
\tilBun_{G,\Del}(\Bl_{x,k}(\Sig)^0)\overset{\tgrq_{\Del}}\to\Bun_{G,\Del}(\Sig)\x \ZZ
$$
are constructed in the obvious
way.

With these notations in order to prove property (i) we must show the following: fix some $a\in\ZZ$ and $\lam\in\Lam$.
Consider the set $\Xi$ of isomorphism classes of all points $\xi\in\tilBun_{G,\Del}(\Bl_{x,k}(\Sig)^0)$ such that

a)$d(\xi)=\lam$.

b) The projection of $\tgrq_{\Del}(\xi)$ to $\ZZ$ is equal to $a$.

Consider now the set $d(\tgrq_{\Del}(\Xi))$ of all possible defects of elements of $\tgrq_{\Del}(\Xi)$. We must show
that this set is finite. To see this, assume that we have a point $\xi$ as above of the form
$(\tcF,\calF_B,\calF_T')$, where $\tcF$ is an element of $\tilBun_{G,\Del}(\Bl_{x,k}(\Sig)^0)$ lying over
some $\calF\in\Bun_{G,\Del}(\Bl_{x,k}(\Sig)^0)$ and
$\calF_B$ and $\calF_T'$ are as before. Let $\calF$ denote the underlying $G$-bundle on $\Bl_{x,k}(\Sig)^0$.
Let $\calG$ denote the unique extension of $\calF|_{\Sig^0}$ to $\Sig$ and let $\calF'$ denote the pull-back
of $\calG$ to $\Bl_{x,k}(\Sig)^0$. Note that $\calF|_{X^0}=\calF'|_{X^0}$, thus $\calF'$ comes with canonical
$\Del$-structure. Moreover, it is clear that $d(\calF',\calF_B,\calF_T')$ is equal to
$d(\grq(\calF,\calF_B,\calF_T'))$. Thus what we have to show is that the defect of $(\calF',\calF_B,\calF_T')$
lies in some finite subset of $\Lam$. This follows from condition b) together with \reft{finite-kapranov}.
\ssec{}{Construction of the Satake isomorphism}We are now in the position to construct the isomorphism claimed in
\reft{main}(2). As was mentioned above the quotient $\tilGam_0\backslash\tilOme$ can be naturally identified with
$\Lam_{\aff}$. Thus the space of $\tilGam_0$-invariant functions on $\tilOme$, supported on finitely many $\tilGam_0$-orbits
is naturally isomorphic to $\CC[\Lam_{\aff}]=\CC(T^{\vee}_{\aff})$. This together with \reft{goodaction} implies that
we the action of $\calH(\tilGam,\Gam_0)$ on $\calF(\tilGam_0\backslash\tilOme)$ gives rise to a homomorphism
$\calH(\tilGam,\tilGam_0)\to \CC(\hatT^{\vee}_{\aff})$. Let $\rho^{\vee}_{\aff}$ denote any element of $\Lam^{\vee}_{\aff}$,
whose scalar product with every simple root of $G_{\aff}^{\vee}$ is equal to 1.
Then Theorem 3.3.5 of \cite{Kap-H} implies that the composition of the above homomorphism  with the shift by
by $q^{-\rho_{\aff}^{\vee}}$ lands in $\CC(\hatT^{\vee}_{\aff})^{W_{\aff}}$. We claim that the resulting homomorphism
\eq{iota}
\iota:\calH(\tilGam^+,\tilGam_0)\to \CC(\hatT^{\vee}_{\aff})^{W_{\aff}}
\eeq
is an isomorphism.

\prop{injective}$\iota$ is injective.
\eprop
\prf
Recall that
$$
\calH(\tilGam^+,\tilGam_0)=\calH(\tilGam^+,\tilGam_0)^{fin}\underset{\CC[v,v^{-1}]}\ten\CC((v))
$$
where
$\calH(\tilGam,\tilGam_0)^{fin}$ is the space of finitely supported functions on
$\tilGam_0\backslash\tilGam^+/\tilGam_0$. Moreover, for any $k>0$ we have
$\tilGam_0\backslash\tilpi^{-1}(k)/\tilGam_0=\Lam_{\aff,k}/W_{\aff}$.
In particular, $\calH(\tilGam,\tilGam_0)^{fin}\simeq \CC[\Lam^+_{\aff}]$
as a vector space.
It is enough to show that the restriction
of $\iota$ to $\calH^{fin}(\tilGam^+,\tilGam_0)$ is injective, since $\iota$ is a morphism of graded
$\CC((v))$-algebras and every graded component of $\calH(\tilGam,\tilGam_0)^{fin}$ is of finite rank
over $\CC[v,v^{-1}]$.
For any $\lam_{\aff}\in\Lam_{\aff}$ let $\del_{\lam_{\aff}}$ denote the characteristic function of
$\tilGam_0s^{\lam_{\aff}}\tilGam_0$.
Then we have to show that the functions $\iota(\del_{\lam})$ are
linearly independent for all $\lam\in\Lam_{\aff}/W_{\aff}$.

To simplify the discussion, let us prove this in the case when $G$ is simply connected (the general
case is similar, but notationally a bit more cumbersome; we shall leave it to the reader).
Recall that when $G$ is simply connected, the quotient $\Lam_{\aff}/W_{\aff}$ can be identified
with the set $\Lam_{\aff}^+$ of dominant coweights of $G_{\aff}$.
Let us define a partial ordering on $\Lam_{\aff}/W_{\aff}$ by saying that $\lam_{\aff}\geq \mu_{\aff}$ if
$\lam_{\aff}-\mu_{\aff}$ is a sum of positive roots of $\grg_{\aff}^{\vee}$.
Then we claim that for any $\lam_{\aff}\in\Lam^+_{\aff}$ we have
\eq{filtration}
\iota(\del_{\lam_{\aff}})=s_{\lam_{\aff}}+\sum\limits_{\mu_{\aff}<\lam_{\aff}} a_{\mu_{\aff}}s_{\mu_{\aff}}.
\end{equation}
It is clear that \refe{filtration} implies that all the $\iota(\del_{\lam_{\aff}})$ are linearly
independent. The proof of \refe{filtration} is a word-by-word
repetition of the corresponding statement for finite-dimensional semi-simple groups (cf. e.g. \cite{Cart}, page 148).
\epr

To prove the surjectivity, let us note that $\iota$ is a morphism of graded $\CC((v))$-algebras. On the other hand, we have
$$
\dim_{\CC((v))}(\calH(\tilGam^+,\tilGam_0)_k)=\#(\Lam_{\aff,k}/W_{\aff})=
\dim_{\CC((v))}(\CC(\hatT^{\vee}_{\aff})^{W_{\aff}}_k).
$$
Thus the injectivity of $\iota$ implies its surjectivity.


\end{document}